\documentstyle[12pt]{article}

\makeatletter

\def\LRA#1#2{\@tempdimb=\c@enumiv\@tempdima%
   \vcenter{\offinterlineskip\halign{##\cr%
   \hfil${\scriptstyle{#1}}$\hfil\crcr%
   \hbox to \@tempdimb{\rightarrowfill}\cr%
   \noalign{\kern-1ex}%
   \hbox to \@tempdimb{\leftarrowfill}\cr%
   \hfil${\scriptstyle{#2}}$\hfil\crcr}}}
\def\RA#1{\@tempdimb=\c@enumiv\@tempdima\vbox{\offinterlineskip%
   \halign{##\cr\hfil${\scriptstyle {#1}}$\hfil\crcr%
   \hbox to \@tempdimb{\rightarrowfill}\cr}}}
\def\LA#1{\@tempdimb=\c@enumiv\@tempdima\vbox{\offinterlineskip%
   \halign{##\cr\hfil${\scriptstyle {#1}}$\hfil\crcr%
   \hbox to \@tempdimb{\leftarrowfill}\cr}}}

\def\diag{\leavevmode\bgroup\setcounter{enumiv}{1}%
   \unitlength1em \@tempdima3em \def\\{\crcr&}\vbox\bgroup%
   \def\multicolumn##1##2{\multispan##1\setcounter{enumiv}{##1}%
   \hfil{##2}\hfil\setcounter{enumiv}{1}}
   \offinterlineskip\halign\bgroup\vrule height.8em depth.7em %
   width0pt##&&\hfil${\displaystyle{##}}$\hfil\cr&}
\def\enddiag{\crcr\egroup\egroup\egroup}
\makeatother
\if@twoside \oddsidemargin 0pt \evensidemargin 0pt \marginparwidth 85pt
\else \oddsidemargin 0pt \evensidemargin 0pt
\fi
\advance\textheight by \topskip
\font\symbolfont=msbm10
\font\teneu=eufm10
\font\egteu=eufm8
\def\dn#1{\mathchoice{\hbox{\teneu #1}}{\hbox{\teneu #1}}%
   {\hbox{\egteu #1}}{\hbox{\egteu #1}}}

\def\nnn{\hbox{\symbolfont N}}
\def\zzz{\hbox{\symbolfont Z}}

\def\depth{\mathop d\nolimits}

\def\cc{\mathop{\cal C}\nolimits}
\def\dd{\mathop{\cal D}\nolimits}

\def\xx{\mathop{\cal X}\nolimits}
\def\yy{\mathop{\cal Y}\nolimits}

\def\pd{\mathop{\rm pd}\nolimits}
\def\height{\mathop{\rm ht}\nolimits}
\def\id{\mathop{\rm id}\nolimits}
\def\gl{\mathop{\rm gl.dim.}\nolimits}

\def\height{\mathop{\rm ht}\nolimits}

\def\spep#1{\mathop{{}^{\bullet}\strut\kern-.1em{#1}}\nolimits}

\def\ext{\mathop{\rm Ext}\nolimits}
\def\hom{\mathop{\rm Hom}\nolimits}
\def\endm{\mathop{\rm End}\nolimits}

\def\add{\mathop{\rm add}\nolimits}
\def\Mod{\mathop{\rm Mod}\nolimits}
\def\soc{\mathop{\rm soc}\nolimits}

\def\ind{\mathop{\rm ind}\nolimits}
\def\cm{\mathop{\rm CM}\nolimits}

\def\tr{\mathop{\rm Tr}\nolimits}

\def\Cok{\mathop{\rm Cok}\nolimits}
\def\Ker{\mathop{\rm Ker}\nolimits}

\def\Im{\mathop{\rm Im}\nolimits}

\def\supp{\mathop{\rm supp}\nolimits}

\def\der#1{\mathop{{\rm R}^{#1}\alpha}}

\def\mod{\mathop{\rm mod}\nolimits}

\def\rank{\mathop{\rm rank}\nolimits}

\def\add{\mathop{\rm add}\nolimits}

\def\cc{\mathop{\cal C}\nolimits}

\def\tor{\mathop{\rm Tor}\nolimits}
\def\Mod{\mathop{\rm Mod}\nolimits}
\def\Cok{\mathop{\rm Cok}\nolimits}
\def\Ker{\mathop{\rm Ker}\nolimits}

\def\depth{\mathop{\rm depth}\nolimits}
\def\pd{\mathop{\rm pd}\nolimits}

\def\id{\mathop{\rm id}\nolimits}
\def\gl{\mathop{\rm gl.dim}\nolimits}

\def\height{\mathop{\rm ht}\nolimits}

\def\height{\mathop{\rm ht}\nolimits}

\def\Gl{\mathop{\rm GL}\nolimits}
\def\ext{\mathop{\rm Ext}\nolimits}

\def\tr{\mathop{\rm Tr}\nolimits}

\def\resdim#1#2{\mathop{{#1}\mbox{-{\rm dim}}\strut\kern2pt {#2}}\nolimits}

\def\XA{1}
\def\XAA{1.1}
\def\XAAA{1.1.1}
\def\XAAB{1.1.2}
\def\XAAC{1.1.3}
\def\XAB{1.2}
\def\XAC{1.3}
\def\XACA{1.3.1}
\def\XAD{1.4}
\def\XADA{1.4.1}
\def\XADB{1.4.2}
\def\XAE{1.5}
\def\XAF{1.6}
\def\XB{2}
\def\XBA{2.1}
\def\XBAA{2.1.1}
\def\XBB{2.2}
\def\XBBA{2.2.1}
\def\XBBB{2.2.2}
\def\XBBC{2.2.3}
\def\XBBD{2.2.4}
\def\XBC{2.3}
\def\XBCA{2.3.1}
\def\XBCB{2.3.2}
\def\XBD{2.4}
\def\XBDA{2.4.1}
\def\XBDB{2.4.2}
\def\XBE{2.5}
\def\XBEA{2.5.1}
\def\XBEB{2.5.2}
\def\XBF{2.6}
\def\XC{3}
\def\XCA{3.1}
\def\XCAA{3.1.1}
\def\XCB{3.2}
\def\XCBA{3.2.1}
\def\XCC{3.3}
\def\XCCA{3.3.1}
\def\XCD{3.4}
\def\XCDA{3.4.1}
\def\XCDB{3.4.2}
\def\XCDC{3.4.3}
\def\XCDD{3.4.4}
\def\XCDE{3.4.5}
\def\XCE{3.5}
\def\XCEA{3.5.1}
\def\XCEB{3.5.2}
\def\XCF{3.6}
\def\XCFA{3.6.1}
\def\XCFB{3.6.2}
\def\XD{4}
\def\XDA{4.1}
\def\XDB{4.2}
\def\XDBA{4.2.1}
\def\XDBB{4.2.2}
\def\XDBC{4.2.3}
\def\XDBD{4.2.4}
\def\XDBE{4.2.5}
\def\XDC{4.3}
\def\XDCA{4.3.1}
\def\XDCB{4.3.2}
\def\XDCC{4.3.3}
\def\XDCD{4.3.4}
\def\XDCE{4.3.5}
\def\XDD{4.4}
\def\XDDA{4.4.1}
\def\XDDB{4.4.2}
\def\XDDC{4.4.3}
\def\XDDD{4.4.4}

\begin{document}
\begin{center}
\vspace*{1cm}{\large Higher dimensional Auslander-Reiten theory on maximal orthogonal subcategories}
\vskip1em Osamu Iyama
\end{center}

{\footnotesize
{\sc Abstract. }We introduce the concept of maximal $(n-1)$-orthogonal subcategories over artin algebras and orders, and develop $(n+1)$-dimensional Auslander-Reiten theory on them. We give the $n$-Auslander-Reiten translation and the $n$-Auslander-Reiten duality, then show the existence of $n$-almost split sequences and $n$-fundamental sequences. We give some examples.}

\vskip.5em
Auslander-Reiten theory, especially the concept of almost split sequences and their existence theorem, is fundamental to study categories which appear in representation theory, for example, modules over artin algebras [ARS][GR][Ri], their functorially finite subcategories [AS][S], their derived categories [H], Cohen-Macaulay modules over Cohen-Macaulay rings [A2,3][Y], lattices over orders [A2,3][RS], and coherent sheaves on projective curves [AR2][GL]. In these Auslander-Reiten theories, the number `2' is quite symbolic. For one thing, almost split sequences give minimal projective resolutions of simple objects of projective dimension `2' in functor categories. For another, Cohen-Macaulay rings and orders of Krull-dimension `2' have fundamental sequences [A4] and provide us with one of the most beautiful situation in representation theory [A4][RV][Y], which is closely related to McKay's observation [Mc] on simple singularities [AV][E]. In this sense, usual Auslander-Reiten theory should be `2-dimensional' theory, and it would be natural to find a setting for higher dimensional Auslander-Reiten theory from the viewpoint of representation theory and non-commutative algebraic geometry (e.g. [V][ArS][GL]).

In this paper, we introduce $(n-1)$-orthogonal subcategories as a natural domain of a higher dimensional Auslander-Reiten theory (\S\XBB) which should be `$(n+1)$-dimensional'. We show that the $n$-Auslander-Reiten translation functor and the $n$-Auslander-Reiten duality can be defined quite naturally for such categories (\S\XBC,\S\XBCA).
Using them, we show that our categories have {\it $n$-almost split sequences} (\S\XCA,\S\XCCA), which are a new generalization of usual almost split sequences and give minimal projective resolutions of simple objects of projective dimension `$n+1$' in functor categories. We also show the existence of {\it $n$-fundamental sequences} for Cohen-Macaulay rings and orders of Krull-dimension `$n+1$' (\S\XCDD). We show that an invariant subring (of Krull-dimension `$n+1$') corresponding to a finite subgroup $G$ of $\Gl(n+1,k)$ has a natural maximal $(n-1)$-orthogonal subcategory (\S\XBE). In the final section, we give a classification of all maximal $1$-orthogonal subcategories for representation-finite selfinjective algebras and representation-finite Gorenstein orders of classical type (\S\XDBB). We show that the number of such subcategories is related to Catalan numbers (\S\XDBC). These results suggest that our higher dimensional Auslander-Reiten theory is far from abstract non-sense.

\vskip-.5em
\begin{center}
{\sc Acknowledgements}
\end{center}

\vskip-.7em
The author thanks Idun Reiten for her kind advice in correcting my English.

\newpage{\bf\XA\ $n$-Auslander-Reiten translation }

In this section, we study subcategories $\xx^{\dn{M}}_{n-1}$ and $\yy^{\dn{M}}_{n-1}$ of $\mod\Lambda$ for an artin algebra or an order $\Lambda$. Using the syzygy functor and the transpose duality, we construct an equivalence $\tau_n:\underline{\xx}^{\dn{M}}_{n-1}\to\overline{\yy}^{\dn{M}}_{n-1}$. This is regarded as an analogue of Auslander-Reiten translation, and plays an important role in the this paper.

Let $\cc$ be an additive category. We denote by $\cc(X,Y)$ the set of morphisms from $X$ to $Y$, and by $fg\in\cc(X,Z)$ the composition of $f\in\cc(X,Y)$ and $g\in\cc(Y,Z)$. We denote by $J_{\cc}$ the Jacobson radical of $\cc$, and $\ind\cc$ the set of isomorphism classes of indecomposable objects in $\cc$. For $X\in\cc$, we denote by $[X]$ the ideal of $\cc$ consisting of morphisms which factor through $X^n$ for some $n\ge0$.

For a noetherian semiperfect ring $\Lambda$, we denote by $J_\Lambda$ the Jacobson radical of $\Lambda$, by $\Mod\Lambda$ the category of left $\Lambda$-modules, and by $\mod\Lambda$ the category of finitely generated left $\Lambda$-modules. We denote by $(\ )^*$ the functor $\hom_\Lambda(\ ,\Lambda):\mod\Lambda\leftrightarrow\mod\Lambda^{op}$. Let $\underline{\mod}\Lambda:=(\mod\Lambda)/[\Lambda]$ and $\underline{\hom}_\Lambda(X,Y)$ the set of morphisms in $\underline{\mod}\Lambda$. For $X\in\mod\Lambda$, take a minimal projective resolution $P_1\stackrel{f_1}{\to}P_0\stackrel{f_0}{\to}X\to0$. Putting $\Omega X:=\Ker f_0$ and $\tr X:=\Cok f_1^*$, we obtain the {\it syzygy functor} $\Omega:\underline{\mod}\Lambda\rightarrow\underline{\mod}\Lambda$ and the {\it transpose duality} $\tr:\underline{\mod}\Lambda\rightarrow\underline{\mod}\Lambda^{op}$ [AB]. For a subcategory $\cc$ of $\mod\Lambda$, we denote by $\underline{\cc}$ the corresponding subcategory of $\underline{\mod}\Lambda$.

\vskip.5em{\bf\XAA\ Definition }
For $X,Y\in\mod\Lambda$, we write $X\perp_n Y$ if $\ext^i_\Lambda(X,Y)=0$ for any $i$ ($0<i\le n$). For full subcategories $\cc$ and $\dd$ of $\mod\Lambda$, we write $\cc\perp_n\dd$ if $X\perp_n Y$ for any $X\in\cc$ and $Y\in\dd$.  Put $\cc^{\perp_n}:=\{X\in\mod\Lambda\ |\ \cc\perp_n X\}$ and ${}^{\perp_n}\cc:=\{X\in\mod\Lambda\ |\ X\perp_n\cc\}$. Put $\xx_n:={}^{\perp_{n}}({}_\Lambda\Lambda)\subseteq\mod\Lambda$, $\xx_n^{op}:={}^{\perp_{n}}(\Lambda_\Lambda)\subseteq\mod\Lambda^{op}$, $\xx_{n,m}:=\xx_n\cap\tr\xx_m^{op}$ and $\xx_{n,m}^{op}:=\xx_n^{op}\cap\tr\xx_m$. One can easily check the following proposition.

\vskip.5em{\bf\XAAA\ Proposition }{\it
(1) $X\in\xx_{n,m}$ if and only if there exists an exact sequence $P_{n+m+1}\stackrel{f_{n+m+1}}{\longrightarrow}\cdots\stackrel{f_0}{\to}P_0$ such that $P_i$ is projective, $P_1^*\stackrel{f_0^*}{\to}\cdots\stackrel{f_{n+m+1}^*}{\longrightarrow}P_{n+m+1}^*$ is exact and $X=\Cok f_m$.

(2) We have the following diagram whose rows are equivalences and columns are dualities:
\[\begin{diag}
\underline{\xx}_{n,0}&\RA{\Omega}&\underline{\xx}_{n-1,1}&\RA{\Omega}&\cdots&\RA{\Omega}&\underline{\xx}_{1,n-1}&\RA{\Omega}&\underline{\xx}_{0,n}\\
\downarrow^{\tr}&&\downarrow^{\tr}&&&&\downarrow^{\tr}&&\downarrow^{\tr}\\
\underline{\xx}_{0,n}^{op}&\LA{\Omega}&\underline{\xx}_{1,n-1}^{op}&\LA{\Omega}&\cdots&\LA{\Omega}&\underline{\xx}_{n-1,1}^{op}&\LA{\Omega}&\underline{\xx}_{n,0}^{op}
\end{diag}\]}

\vskip-1em{\sc Proof }
(1) If $n=0$ or $m=0$, then the assertion holds obviously. Our general assertion is easily reduced to these cases.

(2) Obviously the duality $\tr:\underline{\mod}\Lambda\leftrightarrow\underline{\mod}\Lambda^{op}$ induces a duality $\underline{\xx}_{n,m}\leftrightarrow\underline{\xx}_{m,n}^{op}$. It is easily checked that $\Omega:\underline{\hom}_\Lambda(X,Y)\to\underline{\hom}_\Lambda(\Omega X,\Omega Y)$ is bijective for any $X\in\xx_1$ and $Y\in\mod\Lambda$ (e.g. [A2;7.4]). Thus each $\Omega$ is full faithful, and dense by (1).\rule{5pt}{10pt}

\vskip.5em{\bf\XAAB\ Corollary }{\it For $n,m\ge0$, $\Omega^m\tr\Omega^{n}$ gives a duality $\underline{\xx}_{n,m}\leftrightarrow\underline{\xx}_{n,m}^{op}$ such that $(\Omega^m\tr\Omega^{n})^2$ is isomorphic to the identity functor.}

\vskip.5em{\bf\XAAC\ }Let us prove the following generalization of [A2;3.2].

\vskip.5em{\bf Proposition }{\it
Let $\Lambda$ be a noetherian ring and $n\ge 1$. For any $i$ ($0<i<n$), there exist functorial isomorphisms below for any $X\in\xx_{n-1}$ and $Y\in\mod\Lambda$.\begin{eqnarray*}
\tor_{n-i}^\Lambda(\tr\Omega^{n-1}X,Y)\simeq\ext^i_\Lambda(X,Y),&&\tor_n^\Lambda(\tr\Omega^{n-1}X,Y)\simeq\underline{\hom}_\Lambda(X,Y)
\end{eqnarray*}}

\vskip-1em{\sc Proof }
Let $P_n\stackrel{f_n}{\to}\cdots\stackrel{f_1}{\to}P_0\stackrel{f_0}{\to}X\to0$ be a projective resolution. Then we have a projective resolution $0\to X^*\stackrel{f_0^*}{\to}P_0^*\stackrel{f_1^*}{\to}\cdots\stackrel{f_n^*}{\to}P_n^*\to\tr\Omega^{n-1}X\to0$ of $\tr\Omega^{n-1}X$. 

(i) Since we have a functorial isomorphism $P^*\otimes_\Lambda Y\simeq{}_\Lambda(P,Y)$ for any $Y\in\mod\Lambda$ and projective $P\in\mod\Lambda$, we have the following commutative diagram of complexes.
\[\begin{diag}
\hom_\Lambda(P_{i-1},Y)&\RA{f_i\cdot}&\hom_\Lambda(P_i,Y)&\RA{f_{i+1}\cdot}&\hom_\Lambda(P_{i+1},Y)\\
\parallel&&\parallel&&\parallel\\
P_{i-1}^*\otimes_\Lambda Y&\RA{f_i^*\otimes1}&P_i^*\otimes_\Lambda Y&\RA{f_{i+1}^*\otimes1}&P_{i+1}^*\otimes_\Lambda Y
\end{diag}\]

Comparing the homology of upper and lower sequences, we have a functorial isomorphism $\tor_{n-i}^\Lambda(\tr\Omega^{n-1}X,Y)\simeq\ext^i_\Lambda(X,Y)$.

(ii) We have an exact sequence $X^*\otimes_\Lambda Y\stackrel{a}{\to}\hom_\Lambda(X,Y)\to\underline{\hom}_\Lambda(X,Y)\to0$ of functorial homomorphisms [AB]. Since we have the commutative diagram below of complexes with the exact upper sequence, we obtain a functorial isomorphism $\tor_n^\Lambda(\tr\Omega^{n-1}X,Y)\simeq\underline{\hom}_\Lambda(X,Y)$.
\[\begin{diag} 
0&\RA{}&\hom_\Lambda(X,Y)&\RA{f_0\cdot}&\hom_\Lambda(P_0,Y)&\RA{f_1\cdot}&\hom_\Lambda(P_1,Y)\\
&&\uparrow^a&&\parallel&&\parallel\\
&&X^*\otimes_\Lambda Y&\RA{f_0^*\otimes1}&P_0^*\otimes_\Lambda Y&\RA{f_1^*\otimes1}&P_1^*\otimes_\Lambda Y&\rule{5pt}{10pt}
\end{diag}\]

\vskip-.5em{\bf\XAB\ Definition }
Let $\Lambda$ be an {\it artin algebra} [ARS]. Thus $\Lambda$ contains a central artinian subring $R$ such that $\Lambda$ is a finitely generated $R$-module. Let $I$ be an injective hull of the $R$-module $R/J_R$. Then we have a duality $D:=\hom_R(\ ,I):\mod\Lambda\leftrightarrow\mod\Lambda^{op}$. Let $\overline{\mod}\Lambda:=(\mod\Lambda)/[D\Lambda]$ and $\overline{\hom}_\Lambda(X,Y)$ the set of morphisms in $\overline{\mod}\Lambda$. For a subcategory $\cc$ of $\mod\Lambda$, we denote by $\overline{\cc}$ the corresponding subcategory of $\overline{\mod}\Lambda$. Put $\yy_n:=(D\Lambda)^{\perp_{n}}\subseteq\mod\Lambda$.

\vskip.5em{\bf\XAC\ Definition }
Let $R$ be a complete regular local ring of dimension $d$ and $\Lambda$ an {\it $R$-order} [A2,3][CR]. Thus $\Lambda$ is an $R$-algebra such that $\Lambda$ is a finitely generated projective $R$-module. A typical example of an order is a commutative complete local Cohen-Macaulay ring $\Lambda$ containing a field since such $\Lambda$ contains a complete regular local subring $R$ [M;29.4]. Let $0\to R\to I_0\to\cdots\to I_d\to0$ be a minimal injective resolution of an $R$-module $R$. We denote by $D:=\hom_R(\ ,I_d)$ the Matlis dual.

In this paper, we assume further that $\Lambda$ is an {\it isolated singularity} [A3], namely $\gl\Lambda_\wp=\height\wp$ for any non-maximal prime ideal $\wp$ of $R$. We call a finitely generated left $\Lambda$-module $M$ a {\it Cohen-Macaulay} $\Lambda$-module if it is a projective $R$-module. We denote by $\cm\Lambda$ the category of Cohen-Macaulay $\Lambda$-modules. We have a duality ${D_d}:=\hom_R(\ ,R):\cm\Lambda\leftrightarrow\cm\Lambda^{op}$. Put $\overline{\cm}\Lambda:=(\cm\Lambda)/[{D_d}\Lambda]$. For a subcategory $\cc$ of $\cm\Lambda$, we denote by $\overline{\cc}$ the corresponding subcategory of $\overline{\cm}\Lambda$. Put $\xx^{\cm}_{n}:=\xx_{n}\cap\cm\Lambda$ and $\yy^{\cm}_{n}:=({D_d}\Lambda)^{\perp_{n}}\cap\cm\Lambda$. We collect basic results.

\vskip.5em{\bf\XACA\ }
(1) $\cm\Lambda={}^{\perp_{\infty}}{D_d}\Lambda$ and $\id{}_\Lambda {D_d}\Lambda=d$ [GN1;2.2]. If $0\to X_d\to C_{d-1}\to\cdots\to C_0$ is an exact sequence with $C_i\in\cm\Lambda$, then $X_d\in\cm\Lambda$ [A2;7.2].

(2) $\cm\Lambda=\tr\xx_d^{op}$ [A2;7.9], and $\depth_RX=\min\{i\ge0\ |\ \ext^i_\Lambda(\Lambda/J_\Lambda,X)\neq0\}$ for any $X\in\mod\Lambda$ [GN2;3.2].

(3) $\underline{\hom}_\Lambda(X,Y)$, $\ext^i_\Lambda(X,Y)$ and $\tor_i^\Lambda(Z,X)$ ($i>0$) are finite length $R$-modules for any $X\in\cm\Lambda$, $Y\in\mod\Lambda$ and $Z\in\mod\Lambda^{op}$ [A2;7.6].

\vskip.5em{\sc Proof }
For the convenience of readers, we give an elementary proof of the equality $\cm\Lambda=\tr\xx_d^{op}$ along the argument in [EG;3.8] where the commutative case is treated. The other assertions are more elementary. Fix any $X\in\mod\Lambda$ with $\tr X\in\xx_d^{op}$. Let $\cdots\to P_1\to P_0\to\tr X\to0$ be a projective resolution. Applying $(\ )^*$, we have an exact sequence $0\to X\to P_2^*\to\cdots\to P_{d+1}^*$ with $P_i^*\in\cm\Lambda$. Thus $X\in\cm\Lambda$ by (1).

Fix any $X\in\cm\Lambda$. Since $\Lambda$ is an isolated singularity, $\ext^i_{\Lambda^{op}}(\tr X,\Lambda)$ has finite length for any $i>0$. Let $\cdots\to P_1\to P_0\to\tr X\to0$ be a projective resolution. Since we have an exact sequence $0\to\ext^1_{\Lambda^{op}}(\tr X,\Lambda)\to X$ with $X\in\cm\Lambda$, we have $\ext^1_{\Lambda^{op}}(\tr X,\Lambda)=0$. Now assume $\tr X\perp_{n-1}\Lambda$ for some $n$ ($2\le n\le d$). Then we have an exact sequence $0\to X\to P_2^*\to\cdots\to P_{n-1}^*\to(\Omega^{n}X)^*\to\ext^{n}_{\Lambda^{op}}(\tr X,\Lambda)\to0$ with $\depth_R(\Omega^{n}X)^*\ge2$ and $X,P_i^*\in\cm R$. Considering depth, we obtain $\ext^{n}_{\Lambda^{op}}(\tr X,\Lambda)=0$. Thus we have $\tr X\perp_d\Lambda$ inductively.\rule{5pt}{10pt}

\vskip.5em{\bf\XAD\ Notation }
Throughout this paper, we will treat artin algebras or orders, and use the following notation. Let $\Lambda$ be an artin algebra over $R$ in \XAB, or an order over a complete regular local ring $R$ of dimension $d$ in \XAC. For the artin algebra case, put ${}_\Lambda\dn{M}:=\mod\Lambda$, $\xx^{\dn{M}}_n:=\xx_n$, $\yy^{\dn{M}}_n:=\yy_n$, $d:=0$ and ${D_d}:=D:{}_\Lambda\dn{M}\leftrightarrow{}_{\Lambda^{op}}\dn{M}$. For the order case, put ${}_\Lambda\dn{M}:=\cm\Lambda$, $\xx^{\dn{M}}_n:=\xx^{\cm}_n$ and $\yy^{\dn{M}}_n:=\yy^{\cm}_n$.

In both cases, ${}_\Lambda\dn{M}$ forms a Krull-Schmidt category, namely any object is isomorphic to a finite direct sum of objects whose endomorphism rings are local. We call $X\in{}_\Lambda\dn{M}$ {\it injective} if ${D_d}X$ is a projective $\Lambda^{op}$-module. Then the {\it Nakayama functor} $\nu:={D_d}(\ )^*$ gives an equivalence from projective objects in ${}_\Lambda\dn{M}$ to injective objects in ${}_\Lambda\dn{M}$. The inverse of $\nu$ is given by $\nu^-:=(\ )^*{D_d}$.

\vskip.5em{\bf\XADA\ }
We define the {\it $n$-Auslander-Reiten translation} as in the following theorem, where the usual one is given by $\tau=\tau_1$ and $\tau^-=\tau^-_1$.

\vskip.5em{\bf Theorem }{\it
Let $\Lambda$ be in \XAD\ and $n\ge1$. There exist mutually inverse equivalences $\tau_{n}:={D_d}\Omega^d\tr\Omega^{n-1}:\underline{\xx}^{\dn{M}}_{n-1}\to\overline{\yy}^{\dn{M}}_{n-1}$ and $\tau_{n}^-:=\Omega^d\tr\Omega^{n-1}{D_d}:\overline{\yy}^{\dn{M}}_{n-1}\to\underline{\xx}^{\dn{M}}_{n-1}$. Thus $\tau_n$ gives a bijection from non-projective objects in $\ind\xx^{\dn{M}}_{n-1}$ to non-injective objects in $\ind\yy^{\dn{M}}_{n-1}$, and the inverse is given by $\tau_n^-$.}

\vskip.5em{\sc Proof }
If $\Lambda$ is an artin algebra, then we only have to compose the dualities $\tr\Omega^{n-1}:\underline{\xx}_{n-1}\leftrightarrow\underline{\xx}_{n-1}^{op}$ in \XAAB\ and $D:\underline{\xx}_{n-1}^{op}\leftrightarrow\overline{\yy}_{n-1}$. If $\Lambda$ is an order, then $\xx^{\cm}_{n-1}=\xx_{n-1,d}$ by \XACA(2). Thus $\Omega^d\tr\Omega^{n-1}$ gives a duality $\underline{\xx}^{\cm}_{n-1}\leftrightarrow\underline{\xx}^{\cm}_{n-1}{}^{op}$ by \XAAB, and ${D_d}$ gives a duality ${D_d}:\underline{\xx}^{\cm}_{n-1}{}^{op}\leftrightarrow\overline{\yy}^{\cm}_{n-1}$.\rule{5pt}{10pt}

\vskip.5em{\bf\XADB\ }We note here that the case $d=n+1$ is quite peculiar since these functors $\tau_{n}$ and $\tau_{n}^-$ lift to the Nakayama functors $\nu$ and $\nu^-$ respectively.

\vskip.5em{\bf Theorem }{\it
Let $\Lambda$ be in \XAD\ and $d=n+1\ge2$. Then $\nu={D_d}(\ )^*$ gives an equivalence $\xx^{\dn{M}}_{n-1}\to\yy^{\dn{M}}_{n-1}$ with inverse $\nu^-=(\ )^*{D_d}:\yy^{\dn{M}}_{n-1}\to\xx^{\dn{M}}_{n-1}$, and $\nu$ and $\nu^-$ give lifts of $\tau_n$ and $\tau_n^-$ respectively.}

\vskip.5em{\sc Proof }
Since $(\ )^*=\Omega^2\tr$ holds, the latter assertion follows from the definition of $\tau_n^{\pm}$ and \XAAA. Fix $X\in\xx^{\dn{M}}_{n-1}$. Since $X\in\tr\xx_d^{op}$ by \XACA(2), the exact sequence $0\to\ext^1_\Lambda(\tr X,\Lambda)\to X\to X^{**}\to\ext^2_\Lambda(\tr X,\Lambda)\to0$ [AB] implies that $X$ is reflexive. Similarly, ${D_d}Y$ is reflexive for any $Y\in\yy^{\dn{M}}_{n-1}$, and the former assertion follows.\rule{5pt}{10pt}

\vskip.5em{\bf\XAE\ }We have the following {\it $n$-Auslander-Reiten duality}, where the usual one [A2;8.8] (see also [Y;3.10]) is the case $n=1$.

\vskip.5em{\bf Theorem }{\it
Let $\Lambda$ be in \XAD\ and $n\ge1$. For any $i$ ($0<i<n$), there exist functorial isomorphisms below for any $X\in\xx^{\dn{M}}_{n-1}$, $Y\in\yy^{\dn{M}}_{n-1}$ and $Z\in{}_\Lambda\dn{M}$.
\begin{eqnarray*}
\ext^{n-i}_\Lambda(X,Z)\simeq D\ext^{i}_\Lambda(Z,\tau_nX),&&\underline{\hom}_\Lambda(X,Z)\simeq D\ext^n_\Lambda(Z,\tau_nX)\\
\ext^{n-i}_\Lambda(Z,Y)\simeq D\ext^{i}_\Lambda(\tau_n^-Y,Z),&&\overline{\hom}_\Lambda(Z,Y)\simeq D\ext^n_\Lambda(\tau_n^-Y,Z)
\end{eqnarray*}}

\vskip-1em{\sc Proof }
Notice that we have a functorial isomorphism $\hom_R(\tor_{i}^\Lambda(W,Z),I)=\ext^{i}_\Lambda(W,\hom_R(Z,I))$ ($i\ge0$) for any injective $R$-module $I$ [CE]. Thus, if $\Lambda$ is an artin algebra, then we have functorial isomorphisms $D\tor_{i}^\Lambda(W,Z)=\ext^{i}_\Lambda(W,DZ)=\ext^{i}_\Lambda(Z,DW)$. Putting $W:=\tr\Omega^{n-1}X$, we obtain the assertion by \XAAC. Now let $\Lambda$ be an order. Taking $\hom_R(Z,\ )$, we obtain an exact sequence $0\to {D_d}Z\to{}_R(Z,I_0)\to\cdots\to{}_R(Z,I_d)\to0$. The remark above implies $\ext^{i}_\Lambda(W,\hom_R(Z,I_j))=\hom_R(\tor_{i}^\Lambda(W,Z),I_j)$ ($i,j\ge0$). Since $\tor_{i}^\Lambda(W,Z)$ ($i>0$) has a finite length by \XACA(3), we have $\ext^{i}_\Lambda(W,\hom_R(Z,I_j))=\hom_R(\tor_{i}^\Lambda(W,Z),I_j)=0$ ($i>0$, $j<d$). Thus we have $\hom_R(\tor_{i}^\Lambda(W,Z),I_d)=\ext^{i}_\Lambda(W,\hom_R(Z,I_d))=\ext^{i+d}_\Lambda(W,{D_d}Z)=\ext^{i}_\Lambda(\Omega^dW,{D_d}Z)=\ext^{i}_\Lambda(Z,{D_d}\Omega^dW)$. Putting $W:=\tr\Omega^{n-1}X$, we obtain the assertion by \XAAC.\rule{5pt}{10pt}

\vskip.5em{\bf\XAF\ }The theorem below will be crucial to study sink maps (called minimal right almost split maps in [ARS]) of projective modules.

\vskip.5em{\bf Theorem }{\it
Let $\Lambda$ be in \XAD\ and $X$ a finite length $\Lambda$-module with a projective resolution $\cdots\to P_1\stackrel{f}{\to}P_0\to X\to 0$. Then we have an exact sequence ${}_\Lambda(\ ,\nu P_1)\stackrel{\nu f}{\to} {}_\Lambda(\ ,\nu P_0)\to D\ext^d_\Lambda(X,\ )\to0$ on ${}_\Lambda\dn{M}$.}

\vskip.5em{\sc Proof }
Notice that ${}_\Lambda(\ ,\nu P)={D_d}{}_\Lambda(P,\ )$ holds for any projective $\Lambda$-module $P$. Thus the assertion for $d=0$ follows easily. Assume $d>0$ and fix $Y\in{}_\Lambda\dn{M}$. Then $\ext^d_R(\ext^d_\Lambda(X,Y),R)=D\ext^d_\Lambda(X,Y)$ holds on ${}_\Lambda\dn{M}$. By \XACA(2), we have an exact sequence $0\to{}_\Lambda(P_0,Y)\stackrel{}{\to}\cdots\stackrel{}{\to}{}_\Lambda(P_{d-1},Y)\stackrel{}{\to}{}_\Lambda(\Omega^dX,Y)\stackrel{a}{\to}\ext^d_\Lambda(X,Y)\to0$ on ${}_\Lambda\dn{M}$. If $d=1$, then we have exact sequences $0\to {D_d}{}_\Lambda(\Omega X,Y)\to {D_d}{}_\Lambda(P_0,Y)\to\ext^{1}_R(\ext^1_\Lambda(X,Y),R)\to0$ and ${D_d}{}_\Lambda(P_1,Y)\to {D_d}{}_\Lambda(\Omega X,Y)\to0$, and the assertion follows. Now assume $d\ge2$ and put $Z:=\Ker a$. By an exact sequence $0\to Z\to {}_\Lambda(\Omega^dX,Y)\stackrel{a}{\to}\ext^d_\Lambda(X,Y)\to0$, we have an exact sequence $\ext^{d-1}_R({}_\Lambda(\Omega^dX,Y),R)\to\ext^{d-1}_R(Z,R)\to\ext^{d}_R(\ext^d_\Lambda(X,Y),R)\to\ext^{d}_R({}_\Lambda(\Omega^dX,Y),R)$. Since $\depth_R{}_\Lambda(\Omega^dX,Y)\ge2$ holds, we have $\ext^{i}_R({}_\Lambda(\Omega^dX,Y),R)=0$ ($i\ge d-1$). Thus $\ext^{d-1}_R(Z,R)=\ext^{d}_R(\ext^d_\Lambda(X,Y),R)$. Since $0\to{}_\Lambda(P_0,Y)\stackrel{}{\to}\cdots\stackrel{}{\to}{}_\Lambda(P_{d-1},Y)\to Z\to0$ gives a projective resolution of $R$-modules, we have an exact sequence ${D_d}{}_\Lambda(P_1,Y)\to {D_d}{}_\Lambda(P_1,Y)\to\ext^{d-1}_R(Z,R)\to0$.\rule{5pt}{10pt}

\newpage{\bf\XB\ Maximal $(n-1)$-orthogonal subcategory }

In this section, we keep the notation in \XAD. We introduce maximal $(n-1)$-orthogonal subcategories $\cc$ of ${}_\Lambda\dn{M}$. We show that the $n$-Auslander-Reiten translation defined in \XADA\ induces an equivalence $\underline{\cc}\to\overline{\cc}$, and that $\cc$ has the $n$-Auslander-Reiten duality. The case $d=n+1$ is especially nice, and invariant subrings give such examples.

Let us start by recalling the following useful concept. A morphism $f:X\to Y$ is called {\it right minimal} (resp. {\it left minimal}) if it does not have a direct summand of the form $X\to0$ (resp. $X\to0$) with $X\neq0$ as a complex [ARS].

\vskip.5em{\bf\XBA\ Definition }
Let $\cc$ be a full subcategory of ${}_\Lambda\dn{M}$. We call a morphism $f:C\to X$ a {\it right $\cc$-approximation} of $X$ if $C\in\cc$ and ${}_\Lambda(\ ,C)\stackrel{\cdot f}{\to}{}_\Lambda(\ ,X)\to0$ is exact on $\cc$. We call $\cc$ {\it contravariantly finite} if any $X\in{}_\Lambda\dn{M}$ has a right $\cc$-approximation [AS]. We call a complex ${\bf A}:\cdots\to C_2\stackrel{f_2}{\to}C_1\stackrel{f_1}{\to}C_0\stackrel{f_0}{\to}X$ a {\it right $\cc$-resolution} if $C_i\in\cc$ and $\cdots\to{}_\Lambda(\ ,C_2)\stackrel{\cdot f_2}{\to}{}_\Lambda(\ ,C_1)\stackrel{\cdot f_1}{\to}{}_\Lambda(\ ,C_0)\stackrel{\cdot f_0}{\to}{}_\Lambda(\ ,X)\to0$ is exact on $\cc$. If each $f_i$ is right minimal, then we call ${\bf A}$ {\it minimal}. We write $\resdim{\cc}{X}\le n$ if $X$ has a right $\cc$-resolution with $C_{n+1}=0$. Dually, a {\it left $\cc$-approximation}, a {\it covariantly finite} subcategory, a {\it (minimal) left $\cc$-resolution} and $\resdim{\cc{}^{op}}{X}$ are defined. We call $\cc$ {\it functorially finite} if it is contravariantly and covariantly finite. The following fact immediately follows from Auslander-Buchweitz approximation theory [ABu].

\vskip.5em{\bf\XBAA\ }
Let $\cc$ be a functorially finite subcategory of ${}_\Lambda\dn{M}$. Then any $X\in\mod\Lambda$ has a minimal right (resp. left) $\cc$-resolution, which is unique up to isomorphisms of complexes.

\vskip.5em{\bf\XBB\ Definition }
Let $\cc$ be a functorially finite subcategory of ${}_\Lambda\dn{M}$ and $l\ge0$. We call $\cc$ {\it maximal $l$-orthogonal} if
\[\cc=\cc^{\perp_{l}}\cap{}_\Lambda\dn{M}={}^{\perp_{l}}\cc\cap{}_\Lambda\dn{M}.\]
Then $\cc$ is closed under direct summands, and satisfies $\Lambda\oplus {D_d}\Lambda\in\cc$ and $\cc\subseteq\xx^{\dn{M}}_l\cap\yy^{\dn{M}}_l$. By definition, ${}_\Lambda\dn{M}$ itself is a unique maximal $0$-orthogonal subcategory. If a maximal $l$-orthogonal subcategory $\cc$ is contained in a subcategory $\dd$ of ${}_\Lambda\dn{M}$ satisfying $\dd\perp_l\dd$, then $\cc=\dd$. Moreover, $\cc$ is a maximal $l$-orthogonal subcategory of ${}_\Lambda\dn{M}$ if and only if ${D_d}\cc$ is that of ${}_{\Lambda^{op}}\dn{M}$.

We often use the following easy observation.

\vskip.5em{\bf\XBBA\ }Let $\cc$ be a subcategory of ${}_\Lambda\dn{M}$ and ${\bf A}:0\to X_n\stackrel{f_{n}}{\to}C_{n-1}\stackrel{f_{n-1}}{\to}\cdots\stackrel{f_1}{\to}C_0\stackrel{f_0}{\to}X_0\to0$ an exact sequence with $C_i\in\cc$.

(1) $\ext^1_\Lambda(\ ,X_0)=\ext^{n+1}_\Lambda(\ ,X_n)$ on $\cc^{\perp_{n+1}}$, and $\ext^1_\Lambda(\ ,X_0)\subseteq\ext^{n+1}_\Lambda(\ ,X_n)$ on $\cc^{\perp_n}$.


(2) If $\cc\perp_n\cc$ and ${\bf A}$ is a right $\cc$-resolution, then $X_n\in\cc^{\perp_n}$.

\vskip.5em{\sc Proof }Put $X_i:=\Ker f_{i-1}$. We have an exact sequence $0\to X_i\to C_{i-1}\to X_{i-1}\to0$.

(1) $\ext^1_\Lambda(\ ,X_{0})=\ext^2_\Lambda(\ ,X_{1})=\cdots=\ext^{n+1}_\Lambda(\ ,X_{n})$ holds on ${}^{\perp_{n+1}}\cc$. Similarly, $\ext^1_\Lambda(\ ,X_{0})=\cdots=\ext^n_\Lambda(\ ,X_{n-1})\subseteq\ext^{n+1}_\Lambda(\ ,X_{n})$ holds on ${}^{\perp_{n}}\cc$.


(2) Since $\cc(\ ,C_{i-1})\to\cc(\ ,X_{i-1})\to0$ is exact, $\ext^1_\Lambda(\ ,X_i)=0$ holds on $\cc$ for any $i$ ($1\le i\le n$). Applying (1) to the sequence $0\to X_n\stackrel{f_{n}}{\to}C_{n-1}\stackrel{f_{n-1}}{\to}\cdots\stackrel{}{\to}X_{n-i+1}\to0$, $\ext^i_\Lambda(\ ,X_n)=\ext^1_\Lambda(\ ,X_{n-i+1})=0$ holds on $\cc$ for any $i$ ($1\le i\le n$).\rule{5pt}{10pt}

\vskip.5em{\bf\XBBB\ }We have the following characterizations of maximal $l$-orthogonal subcategories.

\vskip.5em{\bf Proposition }{\it
Let $\cc$ be a functorially finite subcategory of ${}_\Lambda\dn{M}$. Then the conditions (1), (2-$i$) and (3-$i$) are equivalent for any $i$ ($0\le i\le l$).

(1) $\cc$ is maximal $l$-orthogonal.

(2-0) $\resdim{\cc}{X}\le l$ for any $X\in{}_\Lambda\dn{M}$, $\cc\perp_l\cc$ and $\Lambda\oplus D_d\Lambda\in\cc$.

(2-$i$) $\resdim{\cc}{X}\le l-i$ for any $X\in\cc^{\perp_i}\cap{}_\Lambda\dn{M}$, $\cc\perp_l\cc$ and $\Lambda\oplus D_d\Lambda\in\cc$.

(2-$l$) $\cc=\cc^{\perp_l}\cap{}_\Lambda\dn{M}$ and $\Lambda\in\cc$.

(3-0) $\resdim{\cc^{op}}{X}\le l$ for any $X\in{}_\Lambda\dn{M}$, $\cc\perp_l\cc$ and $\Lambda\oplus D_d\Lambda\in\cc$.

(3-$i$) $\resdim{\cc^{op}}{X}\le l-i$ for any $X\in{}^{\perp_i}\cc\cap{}_\Lambda\dn{M}$, $\cc\perp_l\cc$ and $\Lambda\oplus D_d\Lambda\in\cc$.

(3-$l$) $\cc={}^{\perp_l}\cc\cap{}_\Lambda\dn{M}$ and ${D_d}\Lambda\in\cc$.}

\vskip.5em{\sc Proof }
Since (1) is equivalent to (2-$l$)+(3-$l$) by definition, we only have to show that all conditions except (1) are equivalent. We will show that (2-$(i+1)$) implies (2-$i$). For any $X\in\cc^{\perp_i}\cap{}_\Lambda\dn{M}$, take a right $\cc$-resolution $0\to Y\to C_0\to X\to0$, which is exact by $\Lambda\in\cc$. Then $\cc\perp_{l}C_0$ and $\cc\perp_iX$ imply $\cc\perp_{i+1}Y$. Since $\resdim{\cc}{Y}\le l-i-1$ by (2-$(i+1)$), we obtain $\resdim{\cc}{X}\le l-i$. Dually, (3-$(i+1)$) implies (3-$i$).

We will show that (2-0) implies (3-$l$). Fix $X_0\in{}^{\perp_l}\cc\cap{}_\Lambda\dn{M}$. Take a right $\cc$-resolution ${\bf A}:C_{l-1}\stackrel{f_{l-1}}{\to}\cdots\stackrel{f_1}{\to}C_0\stackrel{f_0}{\to}X_0\to0$. Put $X_{l}:=\Ker f_{l-1}$. Since ${\bf A}$ is exact by $\Lambda\in\cc$, we obtain $X_{l}\in\cc^{\perp_l}\cap{}_\Lambda\dn{M}=\cc$ by \XBBA(2) and \XACA(1). 
Applying \XBBA(1) to the sequence $0\to X_l\stackrel{f_{l}}{\to}C_{l-1}\stackrel{f_{l-1}}{\to}\cdots\stackrel{}{\to}X_1\to0$, we have $\ext^1_\Lambda(X_0,X_1)=\ext^{l}_\Lambda(X_0,X_{l})=0$. Thus $f_0$ splits, and we have $X_0\in\cc$. Dually, (3-0) implies (2-$l$).\rule{5pt}{10pt}

\vskip.5em{\bf\XBBC\ }
Especially (2-0) and (3-0) above give the following theorem.

\vskip.5em{\bf Theorem }{\it
Let $\cc$ be a maximal $l$-orthogonal subcategory of ${}_\Lambda\dn{M}$.
Any $X\in{}_\Lambda\dn{M}$ has a minimal right $\cc$-resolution $0\to C_{l}\to\cdots\to C_0\to X\to0$ and a minimal left $\cc$-resolution $0\to X\to C_0^\prime\to\cdots\to C_{l}^\prime\to0$, which are exact.}


\vskip.5em{\bf\XBBD\ }Later we will use the following slightly more general assertion.

\vskip.5em{\bf Lemma }{\it
Let $\cc$ be a maximal $l$-orthogonal subcategory of ${}_\Lambda\dn{M}$.
Assume that $X\in\mod\Lambda$ has an exact sequence $0\to X\to C_{-1}\to\cdots\to C_{-k}$ with $k\ge0$ and $C_i\in{}_\Lambda\dn{M}$. Then $X$ has a minimal right $\cc$-resolution $0\to C_{m}\to\cdots\to C_0\to X\to0$ with $m=\max\{l,d-k\}$.}

\vskip.5em{\sc Proof }By \XBAA, we can take a minimal right $\cc$-resolution $C_{m-1}\stackrel{f_{m-1}}{\to}\cdots\stackrel{f_1}{\to}C_0\to X\to0$, which is exact by $\Lambda\in\cc$. By \XBBA(2) and \XACA(1), $\Ker f_{m-1}\in\cc^{\perp_l}\cap{}_\Lambda\dn{M}=\cc$.\rule{5pt}{10pt}


\vskip.5em{\bf\XBC\ }We have the following {\it $n$-Auslander-Reiten translation} on $\cc$.

\vskip.5em {\bf Theorem }{\it
Let $\cc$ be a maximal $(n-1)$-orthogonal subcategory of ${}_\Lambda\dn{M}$ ($n\ge1$).

(1) $\tau_nX\in\cc$ and $\tau_n^-X\in\cc$ for any $X\in\cc$. 

(2) We have mutually inverse equivalences $\tau_n:\underline{\cc}\to\overline{\cc}$ and $\tau_n^-:\overline{\cc}\to\underline{\cc}$.

(3) $\tau_n$ gives a bijection from non-projective objects in $\ind\cc$ to non-injective objects in $\ind\cc$, and the inverse is given by $\tau_n^-$.}

\vskip.5em{\sc Proof }Recall that $\cc\subseteq\xx^{\dn{M}}_n\cap\yy^{\dn{M}}_n$ holds. By \XAE, for any $X,Y\in\cc$ and $i$ ($0<i<n$), we have $\ext^i_\Lambda(Y,\tau_nX)=D\ext^{n-i}_\Lambda(X,Y)=0$. Thus $\tau_nX\in\cc^{\perp_{n-1}}\cap{}_\Lambda\dn{M}=\cc$ holds. Dually, $\tau_n^-X\in\cc$ holds. The assertion (2) follows from \XADA.\rule{5pt}{10pt}

\vskip.5em{\bf\XBCA\ }We obtain the following {\it $n$-Auslander-Reiten duality} on $\cc$ from \XAE.

\vskip.5em{\bf Theorem }{\it
Let $\cc$ be a maximal $(n-1)$-orthogonal subcategory of ${}_\Lambda\dn{M}$ ($n\ge1$). For any $i$ ($0<i<n$), there exist functorial isomorphisms below for any $X,Y\in\cc$.
\begin{eqnarray*}
\underline{\hom}_\Lambda(X,Y)\simeq D\ext^n_\Lambda(Y,\tau_nX),&&\overline{\hom}_\Lambda(X,Y)\simeq D\ext^n_\Lambda(\tau_n^-Y,X)
\end{eqnarray*}}

\vskip-1em{\bf\XBCB\ }For the case $d=n+1$, \XADB\ implies the following theorem.

\vskip.5em{\bf Theorem }{\it
Assume $d=n+1\ge2$ in \XBC. Then $\nu={D_d}(\ )^*$ gives an equivalence $\cc\to\cc$ with inverse $\nu^-=(\ )^*{D_d}$, and $\nu$ and $\nu^-$ give lifts of $\tau_n$ and $\tau_n^-$ respectively.}

\vskip.5em{\bf\XBD\ Definition }
Let $\cc$ be a Krull-Schmidt category (\XAD). A {\it $\cc$-module} is a contravariant additive functor from $\cc$ to the category of abelian groups. For $\cc$-modules $F$ and $F^\prime$, we denote by $\hom(F,F^\prime)$ the set of natural transformations from $F$ to $F^\prime$. Thus we obtain the abelian category $\Mod\cc$ of $\cc$-modules. The Yoneda embedding $X\mapsto\cc(\ ,X)$ (resp. $\cc(X,\ )$) gives a full faithful functor $\cc\to\Mod\cc$ (resp. $\cc^{op}\to\Mod\cc^{op}$). Moreover, $X\mapsto\cc/J_{\cc}(\ ,X)$ (resp. $\cc/J_{\cc}(X,\ )$) gives a bijection from $\ind\cc$ to isomorphism classes of simple $\cc$-modules (resp. $\cc^{op}$-modules).

We study the socle of $\cc$-modules given by the $n$-th extension groups. The corollaries below will play an important role in the study of $n$-almost split sequences.

\vskip.5em{\bf\XBDA\ Corollary }{\it
Let $\cc$ be a maximal $(n-1)$-orthogonal subcategory of ${}_\Lambda\dn{M}$ ($n\ge1$) and $X\in\ind\cc$.

(1) If $X$ is not projective, then the $\cc^{op}$-module $\ext^n_\Lambda(X,\ )$ has a simple socle with $(\soc_{\cc^{op}}\ext^n_\Lambda(X,\ ))(\tau_nX)\neq0$. If $X$ is not injective, then the $\cc$-module $\ext^n_\Lambda(\ ,X)$ has a simple socle with $(\soc_{\cc}\ext^n_\Lambda(\ ,X))(\tau_n^-X)\neq0$.

(2) $(\soc_{\cc}\ext^n_\Lambda(\ ,\tau_nX))(X)=\soc_{\endm_\Lambda(X)}\ext^n_\Lambda(X,\tau_nX)=\soc_{\endm_\Lambda(\tau_nX)^{op}}\ext^n_\Lambda(X,\tau_nX)=(\soc_{\cc^{op}}\ext^n_\Lambda(X,\ ))(\tau_nX)$.}

\vskip.5em{\sc Proof }
(1) Since $Y:=\tau_n^-X$ is indecomposable non-projective, $\underline{\hom}_\Lambda(Y,\ )$ has a simple top $F:=\cc/J_{\cc}(Y,\ )$. Thus $D\underline{\hom}_\Lambda(Y,\ )=\ext^n_\Lambda(\ ,X)$ has a simple socle $DF$ with $DF(Y)\neq0$.

(2) For any $X\in\cc$, we have an exact functor $\Mod\cc\to\Mod\endm_\Lambda(X)$ given by $F\mapsto F(X)$. For any sub $\endm_\Lambda(X)$-module $M$ of $F(X)$, there exists a sub $\cc$-module $F^\prime$ of $F$ such that $M=F^\prime(X)$. Using this fact, it is easily checked that $(\soc_{\cc}\ext^n_\Lambda(\ ,\tau_nX))(X)$ is a simple socle of the $\endm_\Lambda(X)$-module $\ext^n_\Lambda(X,\tau_nX)$.\rule{5pt}{10pt}

\vskip.5em{\bf\XBDB\ Corollary }{\it
Let $S$ be a simple $\Lambda$-module, $P\to S\to0$ a projective cover and $I:=\nu P$. Then the $\cc^{op}$-module $\ext^d_\Lambda(S,\ )$ has a simple socle with $(\soc_{\cc^{op}}\ext^d_\Lambda(S,\ ))(I)\neq0$. Moreover, $\ext^d_\Lambda(S,I)=(\soc_{\cc^{op}}\ext^d_\Lambda(S,\ ))(I)$ and $\ext^d_\Lambda(S,I^\prime)=0$ for any indecomposable injective $I^\prime\neq I$.}


\vskip.5em{\sc Proof }
Since $D\ext^d_\Lambda(S,\ )$ has a simple top $\cc/J_{\cc}(\ ,I)$ by \XAF, the former assertion follows. Applying the Nakayama functor to \XAF, we have an exact sequence $\cc(\nu^-\ ,P_1)\stackrel{\cdot f}{\to}\cc(\nu^-\ ,P_0)\to D\ext^d_\Lambda(S,\ )\to0$ on injective modules. On the other hand, since $\cc(\ ,P_1)\stackrel{\cdot f}{\to}\cc(\ ,P_0)\to\cc/J_{\cc}(\ ,P_0)\to0$ is exact on projective modules, $D\ext^d_\Lambda(S,\ )=\cc/J_{\cc}(\nu^-\ ,P_0)$ holds on injective modules. Thus the latter assertion follows.\rule{5pt}{10pt}

\vskip.5em{\bf\XBE\ }
Let us give a remarkable example of maximal $(n-1)$-orthogonal subcategories with $d=n+1$. For a field $k$, we call $\sigma\in\Gl(d,k)$ a {\it pseudo reflection} if $\rank(\sigma-1)\le 1$. A finite subgroup $G$ of $\Gl(d,k)$ acts on $\Gamma:=k[[x_1,\cdots,x_d]]$ linearly, and the invariant subring $\Lambda:=\Gamma^G$ contains a complete regular local subring $R$ such that $\Lambda$ is an $R$-order. For the case $d=2$, the following theorem is a well-known result which asserts $\add_\Lambda\Gamma=\cm\Lambda$.

\vskip.5em{\bf Theorem }{\it
Let $k$ be a field of characteristic zero, $G$ a finite subgroup of $\Gl(d,k)$ with $d\ge2$, $\Gamma:=k[[x_1,\cdots,x_d]]$ and $\Lambda:=\Gamma^G$ the invariant subring. Assume that $G$ does not contain any pseudo-reflection except the identity, and that $\Lambda$ is an isolated singularity. Then $\cc:=\add_\Lambda\Gamma$ is a maximal $(d-2)$-orthogonal subcategory of $\cm \Lambda$.}

\vskip.5em{\bf\XBEA\ Proposition }{\it
Let $\Lambda$ be an $R$-order which is an isolated singularity, $X,Y\in\cm\Lambda$ and $2\le n\le d$. Then $\depth_R\hom_\Lambda(X,Y)\ge n$ if and only if $X\perp_{n-2}Y$.}

\vskip.5em{\sc Proof }
Take a maximal number $n$ such that $n\le d$ and $X\perp_{n-2}Y$. Let $0\to\Omega^{n-1}X\to P_{n-2}\to\cdots\to P_0\to X\to 0$ be a projective resolution. Then we have an exact sequence ${\bf A}:0\to\hom_\Lambda(X,Y)\to{}_\Lambda(P_0,Y)\to\cdots\to{}_\Lambda(P_{n-2},Y)\to{}_\Lambda(\Omega^{n-1}X,Y)\to\ext^{n-1}_\Lambda(X,Y)\to0$ with $\depth_R{}_\Lambda(\Omega^{n-1}X,Y)\ge2$ and ${}_\Lambda(P_i,Y)\in\cm R$. If $n=d$, then $\hom_\Lambda(X,Y)\in\cm R$ by ${\bf A}$. Now assume $n<d$. Then $\ext^{n-1}_\Lambda(X,Y)\neq0$. Since $\depth_R\ext^{n-1}_\Lambda(X,Y)=0$ by \XACA(3), we obtain $\depth_R\hom_\Lambda(X,Y)=n$ by ${\bf A}$.\rule{5pt}{10pt}

\vskip.5em{\bf\XBEB\ Proof of \XBE\ }
Since $G$ does not contain any pseudo-reflection except the identity, we have $\endm_\Lambda(\Gamma)=\Gamma*G$ by the argument in [A4] (see also [Y;10.7,10.8]) where the assumption $d=2$ is not used. In particular, $\endm_\Lambda(\Gamma)\in\cm R$ holds. Thus \XBEA\ implies $\Gamma\perp_{d-2}\Gamma$. By \XBBB, we only have to show $\Gamma^{\perp_{d-2}}\cap\cm \Lambda\subseteq\cc$. Fix $X\in \Gamma^{\perp_{d-2}}\cap\cm \Lambda$, and take an injective resolution $0\to X\to I_0\to\cdots\to I_{d-1}$ in $\cm \Lambda$. Then $0\to{}_\Lambda(\Gamma,X)\to{}_\Lambda(\Gamma,I_0)\to\cdots\to{}_\Lambda(\Gamma,I_{d-1})$ is exact. Since the $\Gamma$-module ${D_d}\Gamma$ is free, so is ${}_\Lambda(\Gamma,I_i)$ for any $i$. Since the $\Gamma$-module ${}_\Lambda(\Gamma,X)$ is a $d$-th syzygy, it is free. Thus ${}_\Lambda(\Gamma,X)\in\cc$ holds. Since the injection $\Lambda\to\Gamma$ is a split monomorphism of $\Lambda$-modules [A4] (see also [Y;10.5]), $X={}_\Lambda(\Lambda,X)$ is a direct summand of ${}_\Lambda(\Gamma,X)\in\cc$. Thus $X\in\cc$.\rule{5pt}{10pt}

\vskip.5em{\bf\XBF\ }
Let $\Lambda$ be a {\it selfinjective} artin algebra or a {\it Gorenstein} order [A2][DKR]. Namely we assume that ${D_d}\Lambda$ is a projective $\Lambda$-module. Then ${}_\Lambda\overline{\dn{M}}$ coincides with ${}_\Lambda\underline{\dn{M}}$, and the functors $\tau:{}_\Lambda\underline{\dn{M}}\to{}_\Lambda\underline{\dn{M}}$ and $\Omega:{}_\Lambda\underline{\dn{M}}\to{}_\Lambda\underline{\dn{M}}$ are equivalences such that $\tau\Omega=\Omega\tau$ by \XADA\ and \XAAA. The following easy observation will be used in \S\XD.

\vskip.5em{\bf Proposition }{\it
(1) We have a functorial isomorphism $\underline{\hom}_\Lambda(X,Y)\simeq D\underline{\hom}_\Lambda(Y,\Omega^-\tau X)$ for any $X,Y\in{}_\Lambda\underline{\dn{M}}$.

(2) Let $X\in\ind{}_\Lambda\underline{\dn{M}}$. Then the ${}_\Lambda\underline{\dn{M}}$-module $\underline{\hom}_\Lambda(\ ,X)$ has a simple socle with $(\soc_{{}_\Lambda\underline{\dn{M}}}\underline{\hom}_\Lambda(\ ,X))(\tau^-\Omega X)\neq0$. The ${}_\Lambda\underline{\dn{M}}^{op}$-module $\underline{\hom}_\Lambda(X,\ )$ has a simple socle with $(\soc_{{}_\Lambda\underline{\dn{M}}^{op}}\underline{\hom}_\Lambda(X,\ ))(\Omega^-\tau X)\neq0$.}

\vskip.5em{\sc Proof }
(1) Let $0\to\Omega Y\stackrel{}{\to}P\to Y\to0$ be a projective resolution. Then we have an exact sequence ${}_\Lambda(P,\ )\stackrel{}{\to}{}_\Lambda(\Omega Y,\ )\to\ext^1_\Lambda(Y,\ )\to0$. Since $P$ is injective, we have $\underline{\hom}_\Lambda(\Omega Y,\ )=\ext^1_\Lambda(Y,\ )$. Thus $\underline{\hom}_\Lambda(X,Y)=D\ext^1_\Lambda(Y,\tau X)=D\underline{\hom}_\Lambda(\Omega Y,\tau X)=D\underline{\hom}_\Lambda(Y,\Omega^-\tau X)$ by \XAE.


(2) Immediate from (1).\rule{5pt}{10pt}

\vskip.5em{\bf\XC\ $n$-almost split sequences }

In this section, we keep the notation in \XAD. We introduce $n$-almost split sequences on maximal $(n-1)$-orthogonal subcategories. Our main result is the existence theorem of $n$-almost split sequences, and that of $n$-fundamental sequences for the case $d=n+1$. Then we study the functor category $\mod\cc$, especially we show that $\mod\cc$ has finite global dimension.

\vskip.5em{\bf\XCA\ Definition }
Let $\cc$ be a full subcategory of ${}_\Lambda\dn{M}$.

We call a complex $\cdots\stackrel{f_1}{\to}C_1\stackrel{f_1}{\to}C_0\stackrel{f_0}{\to}X$ with terms in $\cc$ a {\it sink sequence} of $X$ (in $\cc$) if $f_i\in J_{\cc}$ and the following sequence is exact.
\begin{eqnarray*}
\cdots\stackrel{\cdot f_2}{\to}\cc(\ ,C_1)\stackrel{\cdot f_1}{\to}\cc(\ ,C_0)\stackrel{\cdot f_0}{\to}J_{\cc}(\ ,X)\to0
\end{eqnarray*}

We call the above $f_0$ a {\it sink map} of $X$. Dually, we call a complex $X\stackrel{f_0^\prime}{\to}C_0^\prime\stackrel{f_1^\prime}{\to}C_1^\prime\stackrel{g_2^\prime}{\to}\cdots$ with terms in $\cc$ a {\it source sequence} of $X$ (in $\cc$) if $f_i^\prime\in J_{\cc}$ and the following sequence is exact.
\begin{eqnarray*}
\cdots\stackrel{\cdot f_2^\prime}{\to}\cc(C_1^\prime,\ )\stackrel{\cdot f_1^\prime}{\to}\cc(C_0^\prime,\ )\stackrel{\cdot f_0^\prime}{\to}J_{\cc}(X,\ )\to0
\end{eqnarray*}

We call the above $f_0^\prime$ a {\it source map} of $X$ (called {\it minimal left almost split map} in [ARS]). A sink (resp. source) sequence induces a minimal projective resolution of a semisimple $\cc$-module $\cc/J_{\cc}(\ ,X)$ (resp. a semisimple $\cc^{op}$-module $\cc/J_{\cc}(X,\ )$). We call an exact sequence $0\to Y\stackrel{f_n}{\to}C_{n-1}\stackrel{f_{n-1}}{\to}\cdots\stackrel{f_1}{\to}C_0\stackrel{f_0}{\to}X\to0$ (resp. $0\to Y\stackrel{f_n}{\to}C_{n-1}\stackrel{f_{n-1}}{\to}\cdots\stackrel{f_1}{\to}C_0\stackrel{f_0}{\to}X$) with terms in $\cc$ an {\it $n$-almost split sequence} (resp. {\it $n$-fundamental sequence}) if it is a sink sequence of $X$ and a source sequence of $Y$ simultaneously. By definition, a direct sum of $n$-almost split (resp. sink, source) sequences and a direct summand of an $n$-almost split (resp. sink, source) sequences are again $n$-almost split (resp sink, source) sequences.

\vskip.5em{\bf\XCAA\ }
By (1) below, a sink (resp. source) sequence is unique up to isomorphisms of complexes. By (2) below, the study of $n$-almost split sequences is reduced to that of $n$-almost split sequences with indecomposable left and right terms.

\vskip.5em{\bf Proposition }{\it
(1) Let ${\bf A}$ and ${\bf A}^\prime$ be sink (resp. source) sequences of $X$ and $X^\prime$ respectively. Then any isomorphism $X\to X^\prime$ extends to an isomorphism ${\bf A}\to{\bf A}^\prime$.

(2) Let ${\bf A}:0\to Y\stackrel{}{\to}C_{n-1}\stackrel{}{\to}\cdots\stackrel{}{\to}C_0\stackrel{}{\to}X\to0$ and ${\bf A}^\prime:0\to Y^\prime\stackrel{}{\to}C_{n-1}^\prime\stackrel{}{\to}\cdots\stackrel{}{\to}C_0^\prime\stackrel{}{\to}X^\prime\to0$ be $n$-almost split sequences. Then $X\simeq X^\prime$ if and only if $Y\simeq Y^\prime$, and $X$ is indecomposable if and only if $Y$ is indecomposable.}

\vskip.5em{\sc Proof }
(1) Immediate from the definition.

(2) The former assertion follows from (1). Moreover, $X\in\cc$ is indecomposable if and only if $\cc/J_{\cc}(\ ,X)\in\Mod\cc$ is indecomposable if and only if the minimal projective resolution of $\cc/J_{\cc}(\ ,X)$ is indecomposable as a complex if and only if ${\bf A}$ is indecomposable as a complex. Similarly, ${\bf A}$ is indecomposable as a complex if and only if $Y\in\cc$ is indecomposable.\rule{5pt}{10pt}


\vskip.5em{\bf\XCB\ }The next lemma gives the relationship among long exact sequences with terms in $\cc$, projective resolutions of $\cc$-modules and those of $\cc^{op}$-modules.

\vskip.5em{\bf Lemma }{\it
Let $\cc$ be a maximal $(n-1)$-orthogonal subcategory of ${}_\Lambda\dn{M}$ ($n\ge1$) and ${\bf A}:0\to X_n\stackrel{f_n}{\to}C_{n-1}\stackrel{f_{n-1}}{\to}\cdots\stackrel{f_1}{\to}C_0\stackrel{f_0}{\to}X_0\to0$ an exact sequence with terms in $\cc$. Put $X_i:=\Ker f_{i-1}$. Then we have exact sequences
\begin{eqnarray*}
0\to\cc(\ ,X_n)\stackrel{\cdot f_n}{\to}\cc(\ ,C_{n-1})\stackrel{\cdot f_{n-1}}{\to}\cdots\stackrel{\cdot f_1}{\to}\cc(\ ,C_0)\stackrel{\cdot f_0}{\to}\cc(\ ,X_0)\to F\to0\\
0\to\cc(X_0,\ )\stackrel{f_0\cdot}{\to}\cc(C_0,\ )\stackrel{f_1\cdot}{\to}\cdots\stackrel{f_{n-1}\cdot}{\to}\cc(C_{n-1},\ )\stackrel{f_n\cdot}{\to}\cc(X_n,\ )\to G\to0
\end{eqnarray*}
on $\cc$ such that
\begin{eqnarray*}
D(G\circ\tau_n^+)=F\subseteq\ext^{n}_\Lambda(\ ,X_{n}),&&D(F\circ\tau_n^-)=G\subseteq\ext^{n}_\Lambda(X_{0},\ ).
\end{eqnarray*}
If ${\bf A}$ is zero in $\ext^n_\Lambda(X_0,X_n)$, then $f_0$ is a split epimorphism and $f_n$ is a split monomorphism.}

\vskip.5em{\sc Proof }
It follows from \XBBA(1) that $\ext^1_\Lambda(\ ,X_i)=\ext^{n-i+1}_\Lambda(\ ,X_{n})=0$ for any $i$ ($1<i\le n$). This shows that the upper sequence is exact. Again by \XBBA(1), $F=\ext^1_\Lambda(\ ,X_1)=\ext^{n-1}_\Lambda(\ ,X_{n-1})$ holds, and $0\to\ext^{n-1}_\Lambda(\ ,X_{n-1})\to\ext^{n}_\Lambda(\ ,X_{n})\stackrel{\cdot f_n}{\to}\ext^{n}_\Lambda(\ ,C_{n-1})$ is exact on $\cc$. By \XBCA, $\overline{\cc}(C_{n-1},\tau_n^+\ )\stackrel{f_n\cdot}{\to}\overline{\cc}(X_n,\tau_n^+\ )\to DF\to0$ is exact. Thus $\cc(C_{n-1},\ )\stackrel{f_n\cdot}{\to}\cc(X_n,\ )\to D(F\circ\tau_n^-)\to0$ is exact, and we obtain $G=D(F\circ\tau_n^-)$.
The latter assertion is immediate since $1_{X_0}$ is mapped to ${\bf A}$ by the inclusion $F(X_0)\subseteq\ext^{n}_\Lambda(X_0,X_{n})$.\rule{5pt}{10pt}

\vskip.5em{\bf\XCBA\ }
Immediately we have the following proposition, where the case $n=1$ is well-known.

\vskip.5em{\bf Proposition }{\it
Let $\cc$ be a maximal $(n-1)$-orthogonal subcategory of ${}_\Lambda\dn{M}$ ($n\ge1$), $0\to X_n\stackrel{f_n}{\to}C_{n-1}\stackrel{f_{n-1}}{\to}\cdots\stackrel{f_1}{\to}C_0\stackrel{f_0}{\to}X_0\to0$ an exact sequence with terms in $\cc$ and $Y\in\cc$. Then $\cc(Y,C_0)\stackrel{\cdot f_0}{\to}\cc(Y,X_0)\to0$ is exact if and only if $\cc(C_{n-1},\tau_nY)\stackrel{f_n\cdot}{\to}\cc(X_n,\tau_nY)\to0$ is exact.}

\vskip.5em{\bf\XCC\ }
We have the following characterizations of $n$-almost split sequences.

\vskip.5em{\bf Proposition }{\it
Let $\cc$ be a maximal $(n-1)$-orthogonal subcategory of ${}_\Lambda\dn{M}$ ($n\ge1$) and ${\bf A}:0\to Y\stackrel{f_n}{\to}C_{n-1}\stackrel{f_{n-1}}{\to}\cdots\stackrel{f_1}{\to}C_0\stackrel{f_0}{\to}X\to0$ an exact sequence with terms in $\cc$. If $f_i\in J_{\cc}$ and $X,Y\in\ind\cc$, then the conditions (1)--(6) are equivalent.

(1) ${\bf A}$ is an $n$-almost split sequence.

(2) $f_0$ is a sink map.

(3) $f_n$ is a source map.

(4) ${\bf A}$ is in $(\soc_{\cc^{op}}\ext^n_\Lambda(X,\ ))(Y)\backslash\{0\}$.

(5) ${\bf A}$ is in $(\soc_{\cc}\ext^n_\Lambda(\ ,Y))(X)\backslash\{0\}$.

(6) $Y=\tau_nX$ and ${\bf A}$ is in $\soc_{\endm_\Lambda(X)}\ext^n_\Lambda(X,\tau_nX)=\soc_{\endm_\Lambda(\tau_nX)^{op}}\ext^n_\Lambda(X,\tau_nX)$.}

\vskip.5em{\sc Proof }We will borrow the notation in \XCB. By \XCB, ${\bf A}$ is non-zero in $\ext^n_\Lambda(X,Y)$.

(1)$\Leftrightarrow$(2)$\Leftrightarrow$(3) The equalities $F=D(G\circ\tau_n^+)$ and $G=D(F\circ\tau_n^-)$ show that $F$ is simple if and only if $G$ is simple. The assertion follows from the exact sequences in \XCB.

(1)$\Leftrightarrow$(5) By \XCB, we have an inclusion $\delta:F=\Cok(\cdot f_0)\subseteq\ext^n_\Lambda(\ ,Y)$ such that ${\bf A}=\delta(1_X)$. Thus $F$ is a simple $\cc$-module if and only if $F\subseteq\soc_{\cc}\ext^n_\Lambda(\ ,Y)$ if and only if $[{\bf A}]\in(\soc_{\cc}\ext^n_\Lambda(\ ,Y))(X)$ since $F$ is generated by $1_X$.

(5)$\Leftrightarrow$(6) Immediate from \XBDA(2).\rule{5pt}{10pt}

\vskip.5em{\bf\XCCA\ }Now we can prove the main result in this section.

\vskip.5em{\bf Theorem }{\it
Let $\cc$ be a maximal $(n-1)$-orthogonal subcategory of ${}_\Lambda\dn{M}$ ($n\ge1$).

(1) For any non-projective $X\in\ind\cc$, there exists an $n$-almost split sequence $0\to Y\to C_{n-1}\to\cdots\to C_0\to X\to0$ in $\cc$.

(2) For any non-injective $Y\in\ind\cc$, there exists an $n$-almost split sequence $0\to Y\to C_{n-1}\to\cdots\to C_0\to X\to0$ in $\cc$.

(3) Any $n$-almost split sequence $0\to Y\stackrel{}{\to}C_{n-1}\stackrel{}{\to}\cdots\stackrel{}{\to}C_0\stackrel{}{\to}X\to0$ satisfies $Y\simeq\tau_nX$ and $X\simeq\tau_n^-Y$.}

\vskip.5em{\sc Proof }
(1) If $n=1$, then $\cc={}_\Lambda\dn{M}$. By \XBDA(1), we can take a non-split short exact sequence ${\bf A}$ which is in $(\soc_{{}_\Lambda\dn{M}}\ext^1_\Lambda(X,\ ))(\tau_nX)$. Then ${\bf A}$ is a $1$-almost split sequence by \XCC. Now assume $n>1$. We can take a sink map $a:Y\to X$ in ${}_\Lambda\dn{M}$ by the case $n=1$, and a right $\cc$-approximation $b:Z\to Y$ of $Y$. Let $f:C_0\to X$ be a right minimal direct summand of $ba:Z\to X$ as a complex. Then $f$ is a sink map in $\cc$. Since $X$ is non-projective, $f$ is surjective. Thus $X_1:=\Ker f$ is in ${}_\Lambda\dn{M}$. By \XBBC, there exists a minimal right $\cc$-resolution $0\to C_n\stackrel{f_n}{\to}\cdots\stackrel{f_2}{\to}C_1\stackrel{}{\to}X_1\to0$ of $X_1$. By \XCC, the sequence $0\to C_n\stackrel{f_n}{\to}\cdots\stackrel{f_1}{\to}C_0\stackrel{f_0}{\to}X\to0$ is an $n$-almost split sequence.

(3) By \XCC, $(\soc_{\cc^{op}}\ext^n_\Lambda(X,\ ))(Y)\neq0$. By \XBDA(1), $Y\simeq \tau_nX$ and $X\simeq \tau_n^-Y$ hold.\rule{5pt}{10pt}


\vskip.5em{\bf\XCD\ }
To complete the study of minimal projective resolutions of simple $\cc$-modules, we will consider sink sequences of projective modules and source sequences of injective modules. Recall that we denote by $d$ the dimension of $R$ (\XAD).

\vskip.5em{\bf\XCDA\ Theorem }{\it
Let $\cc$ be a maximal $(n-1)$-orthogonal subcategory of ${}_\Lambda\dn{M}$ ($n\ge1$), $X\in\ind\cc$, $F:=\cc/J_{\cc}(\ ,X)$ and $G:=\cc/J_{\cc}(X,\ )$.

(1) $d\le\pd_{\cc}F\le\max\{n,d\}$ if $X$ is projective, and $\pd_{\cc}F=n+1$ otherwise.

(2) $d\le\pd_{\cc^{op}}G\le\max\{n,d\}$ if $X$ is injective, and $\pd_{\cc^{op}}G=n+1$ otherwise.}

\vskip.5em{\sc Proof }
(1)  By \XBBD, there exists a minimal right $\cc$-resolution $0\to C_m\stackrel{f_m}{\to}\cdots\stackrel{f_2}{\to}C_1\to J_\Lambda P\to0$ with $m:=\max\{n,d\}$. Then $0\to C_m\stackrel{f_m}{\to}\cdots\stackrel{f_2}{\to}C_1\to P$ gives a sink sequence of $P$. By considering the depth over $R$, we obtain $C_d\neq0$.

(2) Since a sink sequence ${\bf A}$ of $X\in\cc$ gives a source sequence ${D_d}{\bf A}$ of ${D_d}X\in {D_d}\cc$, the assertion follows from (1).\rule{5pt}{10pt}

\vskip.5em{\bf\XCDB\ }For the case $d\ge n+1$, we can give a more explicit construction than \XCDA.

\vskip.5em{\bf Lemma }{\it
Let $\cc$ be a maximal $(n-1)$-orthogonal subcategory of ${}_\Lambda\dn{M}$ with $d\ge n+1\ge2$. For an injective $I\in\ind\cc$, put $P=P_0:=\nu^-I$ and $S:=P/J_\Lambda P$. Take a minimal projective resolution $0\to\Omega^{d-n}S\to P_{d-n-1}\stackrel{f_{d-n-1}}{\longrightarrow}\cdots\stackrel{f_1}{\to}P_0\to S\to 0$ of $S$ and a minimal right $\cc$-resolution $\cdots\stackrel{f_{d-n+2}}{\longrightarrow}C_{d-n+1}\stackrel{f_{d-n+1}}{\longrightarrow}C_{d-n}\to\Omega^{d-n}S\to0$. Then $C_d\simeq I$ and $C_{d+1}=0$ hold, and the following complex gives a source sequence of $I$.
\[C_d\stackrel{f_d}{\to}\cdots\stackrel{f_{d-n+1}}{\longrightarrow}C_{d-n}\to P_{d-n-1}\stackrel{f_{d-n-1}}{\longrightarrow}\cdots\stackrel{f_1}{\to}P_0\to0\]}

\vskip-1em{\sc Proof }
(i) Put $S_{i}:=\Ker f_{i-1}$. Applying \XBBA(2) to the sequence $0\to S_{d-1}\to C_{d-2}\stackrel{f_{d-2}}{\longrightarrow}\cdots\stackrel{f_{d-n+1}}{\longrightarrow}C_{d-n}\to\Omega^{d-n}S\to0$, we obtain $S_{d-1}\in\cc^{\perp_{n-1}}$. By the dual of \XBBA(1) and \XACA(2), $\ext^{1}_\Lambda(S_{d-1},\ )\subseteq\ext^{n}_\Lambda(S_{d-n},\ )=\ext^{d}_\Lambda(S,\ )$ and $\ext^1_\Lambda(S_i,\ )=\ext^{i+1}_\Lambda(S,\ )=0$ ($0\le i<d-1$) hold on $\cc$. Similarly, $\ext^{i}_\Lambda(S_{d-1},{D_d}\Lambda)=\ext^{i+d-1}_\Lambda(S,{D_d}\Lambda)$ ($i>0$) holds. By \XACA(1), $\ext^{i}_\Lambda(S_{d-1},{D_d}\Lambda)=0$ holds for any $i>1$. By \XBDB, $\ext^{1}_\Lambda(S_{d-1},I^\prime)=0$ holds for any indecomposable injective $I^\prime\neq I$, and $\ext^{1}_\Lambda(S_{d-1},I)$ is a simple $\endm_\Lambda(I)^{op}$-module. Thus we can take a non-split exact sequence ${\bf A}:0\to I\stackrel{a}{\to}X\to S_{d-1}\to0$. By \XBDB\ again, ${\bf A}$ is contained in the simple socle of the $\cc^{op}$-module $\ext^{1}_\Lambda(S_{d-1},\ )$.

(ii) We will show $X\in\cc$.

Since $X\in\cc^{\perp_{n-1}}$ holds by (i), we only have to show $X\in{}_\Lambda\dn{M}$. By \XACA(1) and (i), we only have to show $\ext^{1}_\Lambda(X,I)=0$. We have an exact sequence $\endm_\Lambda(I)\stackrel{\delta}{\to}\ext^{1}_\Lambda(S_{d-1},I){\to}\ext^{1}_\Lambda(X,I)\to0$ of $\endm_\Lambda(I)^{op}$-modules. Since $\ext^{1}_\Lambda(S_{d-1},I)$ is a simple $\Lambda^{op}$-module by (i) and $\delta(1_I)=[{\bf A}]\neq0$, we obtain $\ext^{1}_\Lambda(X,I)=0$.

(iii) Since $\cc\perp_1 I$ holds, ${\bf A}$ is a minimal right $\cc$-resolution of $S_{d-1}$. Thus $C_d\simeq I$ and $C_{d+1}=0$ hold. It remains to show that our sequence gives a source sequence of $I$. We have an exact sequence $\cc(X,\ )\stackrel{a\cdot}{\to}\cc(I,\ )\stackrel{\gamma}{\to}\ext^1_\Lambda(S_{d-1},\ )$ of $\cc^{op}$-modules. The image of $\gamma$ is generated by $\gamma(1_I)=[{\bf A}]$, which is contained in the simple socle of the $\cc^{op}$-module $\ext^{1}_\Lambda(S_{d-1},\ )$. Thus the image of $\gamma$ is simple, and $a$ is the source map of $I$. Since $\ext^1_\Lambda(S_i,\ )=0$ ($0\le i<d-1$) holds on $\cc$ by (i), $0\to\cc(P_0,\ )\stackrel{f_1\cdot}{\to}\cdots\stackrel{f_{d-2}\cdot}{\to}\cc(C_{d-2},\ )\stackrel{}{\to}{}_\Lambda(S_{d-1},\ )\to0$ is exact.\rule{5pt}{10pt}

\vskip.5em{\bf\XCDC\ }
We have the following description of source sequences of injective modules and sink sequences of projective modules for the case $d\ge n+1\ge2$.

\vskip.5em{\bf Theorem }{\it
Let $\cc$ be a maximal $(n-1)$-orthogonal subcategory of ${}_\Lambda\dn{M}$ with $d\ge n+1\ge2$. Take a projective $P\in\ind\cc$ and an injective $I\in\ind\cc$ with $I=\nu P$, and put $S:=P/J_\Lambda P$.

(1) There exists an exact sequence
$0\to C_d\stackrel{}{\to}\cdots\stackrel{}{\to}C_{d-n}\to P_{d-n-1}\stackrel{}{\to}\cdots\stackrel{}{\to}P_0\to S\to0$
with $C_d=I$, $P_0=P$ and $P_i$ projective. This gives a source sequence of $I$.

(2) There exists an exact sequence
$0\to I_0\stackrel{}{\to}\cdots\stackrel{}{\to}I_{d-n-1}\to C'_{d-n}\stackrel{}{\to}\cdots\stackrel{}{\to}C'_d\to S\to0$
with $I_0=I$, $C'_d=P$ and $I_i$ injective. This gives a sink sequence of $P$.}

\vskip.5em{\sc Proof }
(1) follows immediately from \XCDB. We now show (2). Since $D_d\cc$ is a maximal $(n-1)$-orthogonal subcategory of ${}_{\Lambda^{op}}\dn{M}$, we can apply (1) to $D_d\cc$. Thus a source sequence of $D_dP$ in $D_d\cc$ is given by an exact sequence $0\to C'_d\stackrel{}{\to}\cdots\stackrel{}{\to}C'_{d-n}\to P'_{d-n-1}\stackrel{}{\to}\cdots\stackrel{}{\to}P'_0\to S'\to0$ with $C'_d=D_dP$, $P'_0=D_dI$, $P'_i$ projective and $S'$ simple. Since $\ext^i_R(S',R)=0$ for any $i\neq d$, we have an exact sequence $0\to D_dP'_0\stackrel{}{\to}\cdots\stackrel{}{\to}D_dP'_{d-n-1}\to D_dC'_{d-n}\stackrel{}{\to}\cdots\stackrel{}{\to}D_dC'_d\to DS'\to0$ by applying $D_d$. This gives a sink sequence of $P=D_dC'_d$ in $\cc$ with $DS'=S$.\rule{5pt}{10pt}

\vskip.5em{\bf\XCDD\ }The concept of fundamental sequences plays an important role for the case $d=2$ in usual Auslander-Reiten theory [A2;3.6] (see also [Y;11.5]). The following theorem shows that $n$-fundamental sequences exist for the case $d=n+1$.

\vskip.5em{\bf Theorem }{\it
Let $\cc$ be a maximal $(n-1)$-orthogonal subcategory of ${}_\Lambda\dn{M}$ and $d=n+1\ge2$.

(1) For any $X\in\cc$, there exists an $n$-fundamental sequence $0\to Y\stackrel{}{\to}C_{n-1}\stackrel{}{\to}\cdots\stackrel{}{\to}C_0\stackrel{}{\to}X$ with $Y=\nu X$.

(2) For any $Y\in\cc$, there exists an $n$-fundamental sequence $0\to Y\stackrel{}{\to}C_{n-1}\stackrel{}{\to}\cdots\stackrel{}{\to}C_0\stackrel{}{\to}X$ with $X=\nu^-Y$.}

\vskip.5em{\sc Proof }
Recall that $\nu={D_d}(\ )^*$ and $\nu^-=(\ )^*{D_d}$ holds on $\cc$ by \XBCB.

(1) By \XCCA, we only have to consider the case when $X=P$ is indecomposable projective. Since $d=n+1$, the exact sequence $0\to C_d\stackrel{f_d}{\to}\cdots\stackrel{f_2}{\to}C_1\stackrel{f_1}{\to}P_0$ in \XCDB\ gives an $n$-fundamental sequence. 

(2) Since $\nu:\cc\to\cc$ is an equivalence by \XBC, the assertion follows from (1).\rule{5pt}{10pt}

\vskip.5em{\bf\XCDE\ }In [I1], the concept of $\tau$-categories was introduced. It was effectively applied in [I2,3] to derive information on $\mod{}_\Lambda\dn{M}$ from almost split sequences. To apply the theory of $\tau$-categories to ${}_\Lambda\underline{\dn{M}}$ in \S\XD, we need the following observation.

\vskip.5em{\bf Theorem }{\it
${}_\Lambda\underline{\dn{M}}$ and ${}_\Lambda\overline{\dn{M}}$ form $\tau$-categories.}

\vskip.5em{\sc Proof }
For any non-injective $X\in\ind{}_\Lambda\underline{\dn{M}}$, the almost split sequence $0\to X\to\theta^-X\to\tau^-X\to0$ induces an exact sequence ${}_\Lambda\underline{\dn{M}}(\tau^-X,\ )\to{}_\Lambda\underline{\dn{M}}(\theta^-X,\ )\to J_{{}_\Lambda\underline{\dn{M}}}(X,\ )\to0$ by [I2;1.3]. If $\tau^-X\in\ind{}_\Lambda\underline{\dn{M}}$, then ${}_\Lambda\underline{\dn{M}}(\ ,X)\to{}_\Lambda\underline{\dn{M}}(\ ,\theta^-X)\to J_{{}_\Lambda\underline{\dn{M}}}(\ ,\tau^-X)\to0$ is exact again by [I2;1.3]. Fix injective $I\in\ind{}_\Lambda\underline{\dn{M}}$ and take the source sequence $I\stackrel{f_d}{\to}C_{d-1}\to P_{d-2}\stackrel{f_{d-2}}{\longrightarrow}\cdots\stackrel{f_1}{\to}P_0$ of $I$ constructed in \XCDB\ for $n=1$. Since $P_{d-2}$ is projective, ${}_\Lambda\underline{\dn{M}}(C_{d-1},\ )\stackrel{f_d\cdot}{\longrightarrow}J_{{}_\Lambda\underline{\dn{M}}}(I,\ )$ is an isomorphism by [I2;1.3]. Thus ${}_\Lambda\underline{\dn{M}}$ forms a $\tau$-category, and the equivalence $\tau:{}_\Lambda\underline{\dn{M}}\to{}_\Lambda\overline{\dn{M}}$ shows that so is ${}_\Lambda\overline{\dn{M}}$.\rule{5pt}{10pt}

\vskip.5em{\bf\XCE\ Definition }
Let $\cc$ be an additive category. We say that a $\cc$-module $F$ is {\it finitely presented} if there exists an exact sequence $\cc(\ ,Y)\stackrel{\cdot f}{\to}\cc(\ ,X)\to F\to0$. We denote by $\mod\cc$ the category of finitely presented $\cc$-modules. The Yoneda embedding $X\mapsto\cc(\ ,X)$ (resp. $X\mapsto\cc(X,\ )$) gives an equivalence from $\cc$ (resp. $\cc^{op}$) to the category of finitely presented projective $\cc$-modules (resp. $\cc^{op}$-modules).

If $\cc$ is an $R$-category for $R$ in \XAD, then $F(X)$ has an $R$-module structure naturally for any $\cc$-module $F$ and $X\in\cc$. Thus we have a functor $D:\Mod\cc\leftrightarrow\Mod\cc^{op}$ by composing with $D:\Mod R\to\Mod R$ in \XAD. We call $\cc$ a {\it dualizing $R$-variety} if $D$ induces a duality $\mod\cc\leftrightarrow\mod\cc^{op}$ [AR1]. In this case, $\mod\cc$ and $\mod\cc^{op}$ are closed under kernels since they are always closed under cokernels [A1;2.1]. Thus $\mod\cc$ (resp. $\mod\cc^{op}$) is an abelian subcategory of $\Mod\cc$ (resp. $\Mod\cc^{op}$).

The following proposition strengthen \XCB.

\vskip.5em{\bf\XCEA\ Proposition }{\it
Let $\cc$ be a maximal $(n-1)$-orthogonal subcategory of ${}_\Lambda\dn{M}$ ($n\ge1$). For any $F\in\mod\underline{\cc}$, there exists a unique exact sequence $0\to C_{n+1}\stackrel{f_{n+1}}{\to}C_{n}\stackrel{f_{n}}{\to}\cdots\stackrel{f_2}{\to}C_1\stackrel{f_1}{\to}C_0\to0$ such that $C_i\in\cc$, $f_i\in J_{\cc}$ and the following sequences are minimal projective resolutions.
\begin{eqnarray*}
&0\to\cc(\ ,C_{n+1})\stackrel{\cdot f_{n+1}}{\to}\cc(\ ,C_{n})\stackrel{\cdot f_{n}}{\to}\cdots\stackrel{\cdot f_2}{\to}\cc(\ ,C_1)\stackrel{\cdot f_1}{\to}\cc(\ ,C_0)\to F\to0&\\
&0\to\cc(C_0,\ )\stackrel{f_1\cdot}{\to}\cc(C_1,\ )\stackrel{f_2\cdot}{\to}\cdots\stackrel{f_{n}\cdot}{\to}\cc(C_{n},\ )\stackrel{f_{n+1}\cdot}{\to}\cc(C_{n+1},\ )\to D(F\circ\tau_n^-)\to0&
\end{eqnarray*}}

\vskip-1em{\sc Proof }
Take a minimal projective resolution $\cc(\ ,C_1)\stackrel{\cdot f_1}{\to}\cc(\ ,C_0)\to F\to0$. Then $F(\Lambda)=0$ implies that $f_1$ is surjective. Thus $X_2:=\Ker f_1$ is in ${}_\Lambda\dn{M}$. By \XBBC, there exists a minimal right $\cc$-resolution $0\to C_{n+1}\stackrel{f_{n+1}}{\to}\cdots\stackrel{f_3}{\to}C_2\to X_2\to0$. Then the exact sequence $0\to C_{n+1}\stackrel{f_{n+1}}{\to}\cdots\stackrel{f_1}{\to}C_0\to0$ gives a minimal projective resolution $0\to\cc(\ ,C_{n+1})\stackrel{\cdot f_{n+1}}{\to}\cdots\stackrel{\cdot f_1}{\to}\cc(\ ,C_0)\to F\to0$. By \XCB, the lower sequence is exact.\rule{5pt}{10pt}

\vskip.5em{\bf\XCEB\ }We have the following description of injective $\underline{\cc}$-modules.

\vskip.5em{\bf Theorem }{\it
Let $\cc$ be a maximal $(n-1)$-orthogonal subcategory of ${}_\Lambda\dn{M}$ ($n\ge1$). 

(1) $\underline{\cc}$ and $\overline{\cc}$ are dualizing $R$-varieties.

(2) $X\mapsto\ext^n_\Lambda(\ ,X)$ gives an equivalence from $\overline{\cc}$ to the category of finitely presented injective $\underline{\cc}$-modules, and $X\mapsto\ext^n_\Lambda(X,\ )$ gives an equivalence from $\underline{\cc}^{op}$ to the category of finitely presented injective $\overline{\cc}^{op}$-modules.}

\vskip.5em{\sc Proof }
(1) If $F\in\mod\underline{\cc}$, then $DF\in\mod\underline{\cc}^{op}$ by \XCEA. By \XACA(3), $D:\mod\underline{\cc}\leftrightarrow\mod\underline{\cc}^{op}$ is in fact a duality.

(2) We only show the former assertion. Since $\ext^n_\Lambda(\ ,X)=D\underline{\hom}_\Lambda(\tau_n^-X,\ )$ holds by \XBCA, our functor $X\mapsto\ext^n_\Lambda(\ ,X)$ is a composition of the equivalence $\tau_n^-:\overline{\cc}\to\underline{\cc}$, the duality $X\mapsto\underline{\cc}(X,\ )$ from $\underline{\cc}$ to the category of finitely presented projective $\underline{\cc}^{op}$-modules, and the duality $D:\mod\underline{\cc}^{op}\to\mod\underline{\cc}$.\rule{5pt}{10pt}

\vskip.5em{\bf\XCF\ Definition }
Let $\cc$ be a functorially finite subcategory of ${}_\Lambda\dn{M}$. Since $\cc$ has pseudokernels by Auslander-Buchweitz approximation theory [ABu][AS], $\mod\cc$ (resp. $\mod\cc^{op}$) is an abelian subcategory of $\Mod\cc$ (resp. $\Mod\cc^{op}$) [A1;2.1]. For $F\in\mod\cc$, take a projective resolution $\cc(\ ,Y)\stackrel{\cdot f}{\to}\cc(\ ,X)\to F\to0$. Then $\alpha F$ is defined by the exact sequence $0\to\alpha F\to\cc(X,\ )\stackrel{f\cdot}{\to}\cc(Y,\ )$. It is easily shown that $\alpha$ gives a left exact functor $\alpha:\mod\cc\rightarrow\mod\cc^{op}$, and we define $\alpha:\mod\cc^{op}\rightarrow\mod\cc$ dually. We denote by $\der{n}:\mod\cc\leftrightarrow\mod\cc^{op}$ their $n$-th derived functor [FGR].

\vskip.5em{\bf\XCFA\ }
Let us calculate $\der{i}F$ for $\underline{\cc}$-modules $F$.

\vskip.5em{\bf Theorem }{\it
Let $\cc$ be a maximal $(n-1)$-orthogonal subcategory of ${}_\Lambda\dn{M}$. 

(1) Any $F\in\mod\underline{\cc}$ satisfies $\pd_{\cc}F\le n+1$ and $\der{i}F=0$ for any $i\neq n+1$.

(2) $\der{n+1}$ gives a duality $\mod\underline{\cc}\leftrightarrow\mod\overline{\cc}^{op}$, and the equivalence $D\der{n+1}:\mod\underline{\cc}\leftrightarrow\mod\overline{\cc}$ coincides with the equivalence induced by $\tau_n:\underline{\cc}\to\overline{\cc}$.}

\vskip.5em{\sc Proof }
(1) and the latter assertion of (2) follows from \XCEA. The former assertion of (2) follows from (1) since $\der{n+1}F=\tr\Omega^nF$ holds (cf. \XAAA(2)).\rule{5pt}{10pt}

\vskip.5em{\bf\XCFB\ }
We now calculate the global dimension of functor categories.

\vskip.5em{\bf Theorem }{\it
Let $\Lambda$ be a maximal $(n-1)$-orthogonal subcategory of ${}_\Lambda\dn{M}$ ($n\ge1$). 

(1) If $\gl\Lambda\neq d$, then $\gl(\mod\cc)=\gl(\mod\cc^{op})=\max\{n+1,d\}$. If $\gl\Lambda=d$, then $\gl(\mod\cc)=\gl(\mod\cc^{op})=d$.

(2) If $F\in\mod\cc$ (resp. $\mod\cc^{op}$) is simple with $\pd F=n+1$, then $\der{i}F=0$ for any $i\neq n+1$ and $\der{n+1}F$ is simple with $\pd\der{n+1}F=n+1$.}

\vskip.5em{\sc Proof }
(1) Take $F\in\mod\cc$ and a projective resolution $\cc(\ ,C_1)\stackrel{f_1}{\to}\cc(\ ,C_0)\to F\to0$. Put $X_2:=\Ker f_1$. By \XBBD, there exists a minimal right $\cc$-resolution $0\to C_{m}\stackrel{f_{m}}{\to}\cdots\stackrel{f_3}{\to}C_2\to X_2\to0$ with $m:=\max\{n+1,d\}$. Thus $\pd F\le m$, and we obtain $\gl(\mod\cc)\le m$. If $\gl\Lambda=d$, then ${}_\Lambda\dn{M}$ consists of projective $\Lambda$-modules [A3;1.5], and the assertion follows. Assume $\gl\Lambda\neq d$. Then there exists non-projective $X\in\ind\cc$, and $\pd_{\cc}\cc/J_{\cc}(\ ,X)=n+1$ by \XCDA(1). If $d>n+1$, then any projective $\Lambda$-module $P$ satisfies $\pd_{\cc}\cc/J_{\cc}(\ ,P)=d$ by \XCDA(1).

(2) Put $F=\cc/J_{\cc}(\ ,X)$ for $X\in\ind\cc$. If $X$ is not projective, then the existence of $n$-almost split sequence implies the assertion. If $X$ is projective, then \XCDA(1) implies $d=n+1$. Thus the existence of $n$-fundamental sequences implies the assertion.\rule{5pt}{10pt}

\vskip.5em{\bf\XD\ Selfinjective algebras and Gorenstein orders }

Throughout this section, we keep the notation in \XAD, and we assume that $\Lambda$ is a {\it representation-finite selfinjective} artin algebra or a {\it representation-finite Gorenstein} order (\XBF). In other words, we assume that ${D_d}\Lambda$ is a projective $\Lambda$-module and $\ind{}_\Lambda\dn{M}$ is a finite set. We will classify all maximal $1$-orthogonal subcategories of ${}_\Lambda\dn{M}$. We denote by $\dn{A}({}_\Lambda\underline{\dn{M}})$ the stable Auslander-Reiten quiver of $\Lambda$. In [R2], Riedtmann introduced the quiver $\zzz\Delta$ associated to a graph $\Delta$. The well-known theorem [R2][HPR] below also holds for our $\Lambda$ since the function $\rank_R$ on ${}_\Lambda\dn{M}$ gives a subadditive function on $\dn{A}({}_\Lambda\underline{\dn{M}})$.

\vskip.5em{\bf\XDA\ Proposition }{\it $\dn{A}({}_\Lambda\underline{\dn{M}})$ is a finite disjoint union $\coprod_{k}(\zzz\Delta_k/G_k)$ with Dynkin diagrams $\Delta_k$ and automorphism groups $G_k$ of $\zzz\Delta_k$.}

\vskip.5em
We denote by $(\zzz\Delta)_0$ (resp. $(\zzz\Delta/G)_0$) the set of vertices of the quiver $\zzz\Delta$ (resp. $\zzz\Delta/G$). Our goal in this section is to determine configurations of maximal $1$-orthogonal subcategories of ${}_\Lambda\dn{M}$ as subsets of $\coprod_k(\zzz\Delta_k/G_k)_0$. For a simplicity reason, {\it we assume that each $\Delta_k$ is a classical Dynkin diagram $A_m$, $B_m$, $C_m$ or $D_{m+1}$} though our method explained in \XDD\ works for arbitrary Dynkin diagram $\Delta$.

\vskip.5em{\bf\XDB\ Definition }
We introduce a coordinate on $\zzz\Delta$ by the diagram below, where we denote $\stackrel{(2,1)}{\longrightarrow}$ or $\stackrel{(1,2)}{\longrightarrow}$ by thick arrows for the case $\Delta=B_m$ or $C_m$.
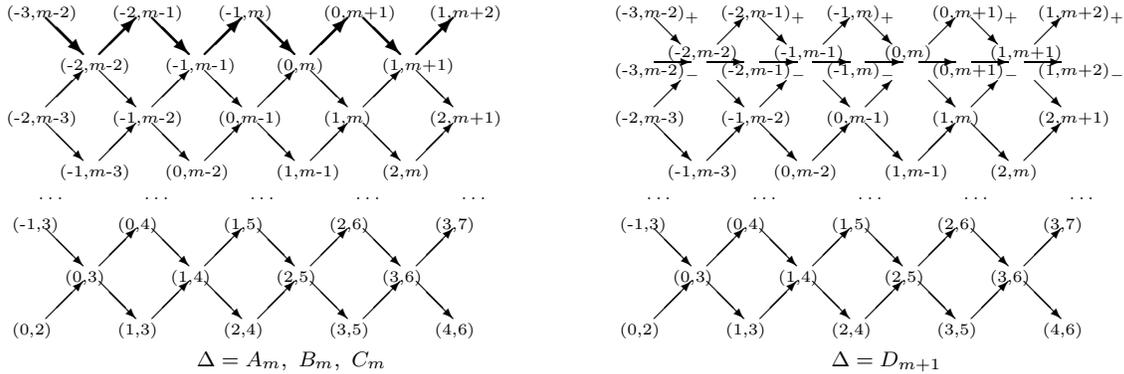
\[\begin{picture}(230,120)
\put(60,-15){\scriptsize$\Delta=A_m,\ B_m,\ C_m$}
\put(3,3){\vector(1,1){14}}
\put(3,37){\vector(1,-1){14}}
\put(3,77){\vector(1,-1){14}}
\put(3,83){\vector(1,1){14}}\thicklines
\put(3,117){\vector(1,-1){14}}\thinlines
\put(-10,-3){\tiny(0,2)}
\put(-10,37){\tiny(-1,3)}
\put(-12,77){\tiny(-2,$m$-3)}
\put(-12,117){\tiny(-3,$m$-2)}

\put(23,17){\vector(1,-1){14}}
\put(23,23){\vector(1,1){14}}
\put(23,63){\vector(1,1){14}}
\put(23,97){\vector(1,-1){14}}\thicklines
\put(23,103){\vector(1,1){14}}\thinlines
\put(10,17){\tiny(0,3)}
\put(8,57){\tiny(-1,$m$-3)}
\put(8,97){\tiny(-2,$m$-2)}

\put(43,3){\vector(1,1){14}}
\put(43,37){\vector(1,-1){14}}
\put(43,77){\vector(1,-1){14}}
\put(43,83){\vector(1,1){14}}\thicklines
\put(43,117){\vector(1,-1){14}}\thinlines
\put(30,-3){\tiny(1,3)}
\put(30,37){\tiny(0,4)}
\put(28,77){\tiny(-1,$m$-2)}
\put(28,117){\tiny(-2,$m$-1)}

\put(63,17){\vector(1,-1){14}}
\put(63,23){\vector(1,1){14}}
\put(63,63){\vector(1,1){14}}
\put(63,97){\vector(1,-1){14}}\thicklines
\put(63,103){\vector(1,1){14}}\thinlines
\put(50,17){\tiny(1,4)}
\put(48,57){\tiny(0,$m$-2)}
\put(48,97){\tiny(-1,$m$-1)}

\put(83,3){\vector(1,1){14}}
\put(83,37){\vector(1,-1){14}}
\put(83,77){\vector(1,-1){14}}
\put(83,83){\vector(1,1){14}}\thicklines
\put(83,117){\vector(1,-1){14}}\thinlines
\put(70,-3){\tiny(2,4)}
\put(70,37){\tiny(1,5)}
\put(68,77){\tiny(0,$m$-1)}
\put(68,117){\tiny(-1,$m$)}

\put(103,17){\vector(1,-1){14}}
\put(103,23){\vector(1,1){14}}
\put(103,63){\vector(1,1){14}}
\put(103,97){\vector(1,-1){14}}\thicklines
\put(103,103){\vector(1,1){14}}\thinlines
\put(90,17){\tiny(2,5)}
\put(90,57){\tiny(1,$m$-1)}
\put(90,97){\tiny(0,$m$)}

\put(123,3){\vector(1,1){14}}
\put(123,37){\vector(1,-1){14}}
\put(123,77){\vector(1,-1){14}}
\put(123,83){\vector(1,1){14}}\thicklines
\put(123,117){\vector(1,-1){14}}\thinlines
\put(110,-3){\tiny(3,5)}
\put(110,37){\tiny(2,6)}
\put(108,77){\tiny(1,$m$)}
\put(108,117){\tiny(0,$m$+1)}

\put(143,17){\vector(1,-1){14}}
\put(143,23){\vector(1,1){14}}
\put(143,63){\vector(1,1){14}}
\put(143,97){\vector(1,-1){14}}\thicklines
\put(143,103){\vector(1,1){14}}\thinlines
\put(130,17){\tiny(3,6)}
\put(130,57){\tiny(2,$m$)}
\put(130,97){\tiny(1,$m$+1)}

\put(150,-3){\tiny(4,6)}
\put(150,37){\tiny(3,7)}
\put(148,77){\tiny(2,$m$+1)}
\put(148,117){\tiny(1,$m$+2)}

\put(0,47){\tiny$\cdots$}
\put(40,47){\tiny$\cdots$}
\put(80,47){\tiny$\cdots$}
\put(120,47){\tiny$\cdots$}
\put(160,47){\tiny$\cdots$}
\end{picture}
\begin{picture}(180,120)
\put(70,-15){\scriptsize$\Delta=D_{m+1}$}
\put(3,3){\vector(1,1){14}}
\put(3,37){\vector(1,-1){14}}
\put(3,77){\vector(1,-1){14}}
\put(3,83){\vector(1,1){10}}
\put(3,117){\vector(1,-1){10}}
\put(3,100){\vector(1,0){14}}
\put(-10,-3){\tiny(0,2)}
\put(-10,37){\tiny(-1,3)}
\put(-12,77){\tiny(-2,$m$-3)}
\put(-12,117){\tiny(-3,$m$-2)$_+$}
\put(-12,95){\tiny(-3,$m$-2)$_-$}

\put(23,17){\vector(1,-1){14}}
\put(23,23){\vector(1,1){14}}
\put(23,63){\vector(1,1){14}}
\put(27,93){\vector(1,-1){10}}
\put(27,107){\vector(1,1){10}}
\put(23,100){\vector(1,0){14}}
\put(10,17){\tiny(0,3)}
\put(8,57){\tiny(-1,$m$-3)}
\put(8,102){\tiny(-2,$m$-2)}

\put(43,3){\vector(1,1){14}}
\put(43,37){\vector(1,-1){14}}
\put(43,77){\vector(1,-1){14}}
\put(43,83){\vector(1,1){10}}
\put(43,117){\vector(1,-1){10}}
\put(43,100){\vector(1,0){14}}
\put(30,-3){\tiny(1,3)}
\put(30,37){\tiny(0,4)}
\put(28,77){\tiny(-1,$m$-2)}
\put(28,117){\tiny(-2,$m$-1)$_+$}
\put(28,95){\tiny(-2,$m$-1)$_-$}

\put(63,17){\vector(1,-1){14}}
\put(63,23){\vector(1,1){14}}
\put(63,63){\vector(1,1){14}}
\put(67,93){\vector(1,-1){10}}
\put(67,107){\vector(1,1){10}}
\put(63,100){\vector(1,0){14}}
\put(50,17){\tiny(1,4)}
\put(48,57){\tiny(0,$m$-2)}
\put(48,102){\tiny(-1,$m$-1)}

\put(83,3){\vector(1,1){14}}
\put(83,37){\vector(1,-1){14}}
\put(83,77){\vector(1,-1){14}}
\put(83,83){\vector(1,1){10}}
\put(83,117){\vector(1,-1){10}}
\put(83,100){\vector(1,0){14}}
\put(70,-3){\tiny(2,4)}
\put(70,37){\tiny(1,5)}
\put(68,77){\tiny(0,$m$-1)}
\put(68,117){\tiny(-1,$m$)$_+$}
\put(68,95){\tiny(-1,$m$)$_-$}

\put(103,17){\vector(1,-1){14}}
\put(103,23){\vector(1,1){14}}
\put(103,63){\vector(1,1){14}}
\put(107,93){\vector(1,-1){10}}
\put(107,107){\vector(1,1){10}}
\put(103,100){\vector(1,0){14}}
\put(90,17){\tiny(2,5)}
\put(90,57){\tiny(1,$m$-1)}
\put(90,102){\tiny(0,$m$)}

\put(123,3){\vector(1,1){14}}
\put(123,37){\vector(1,-1){14}}
\put(123,77){\vector(1,-1){14}}
\put(123,83){\vector(1,1){10}}
\put(123,117){\vector(1,-1){10}}
\put(123,100){\vector(1,0){14}}
\put(110,-3){\tiny(3,5)}
\put(110,37){\tiny(2,6)}
\put(108,77){\tiny(1,$m$)}
\put(108,117){\tiny(0,$m$+1)$_+$}
\put(108,95){\tiny(0,$m$+1)$_-$}

\put(143,17){\vector(1,-1){14}}
\put(143,23){\vector(1,1){14}}
\put(143,63){\vector(1,1){14}}
\put(147,93){\vector(1,-1){10}}
\put(147,107){\vector(1,1){10}}
\put(143,100){\vector(1,0){14}}
\put(130,17){\tiny(3,6)}
\put(130,57){\tiny(2,$m$)}
\put(130,102){\tiny(1,$m$+1)}

\put(150,-3){\tiny(4,6)}
\put(150,37){\tiny(3,7)}
\put(148,77){\tiny(2,$m$+1)}
\put(148,117){\tiny(1,$m$+2)$_+$}
\put(148,95){\tiny(1,$m$+2)$_-$}

\put(0,47){\tiny$\cdots$}
\put(40,47){\tiny$\cdots$}
\put(80,47){\tiny$\cdots$}
\put(120,47){\tiny$\cdots$}
\put(160,47){\tiny$\cdots$}
\end{picture}\]

\vskip.5em
In each case, the translation $\tau$ is defined by $\tau(i,j)=(i-1,j-1)$ and $\tau(i,i+m+1)_{\pm}=(i-1,i+m)_{\pm}$. Define an automorphism $\sigma$ of $\zzz\Delta$ as follows. If $\Delta=A_m$, $B_m$ or $C_m$, then $\sigma:=1_{\zzz\Delta}$. If $\Delta=D_{m+1}$, then $\sigma(i,i+m+1)_{\pm}:=(i,i+m+1)_{\mp}$ and other points are fixed.
For $x\in(\zzz\Delta)_0$, we will define subsets $H^-(x)$ and $H^+(x)$ of $(\zzz\Delta)_0$. If $\Delta=A_m$, then $H^{\pm}(x)$ consists of points in the rectangles (and their boundary) below.
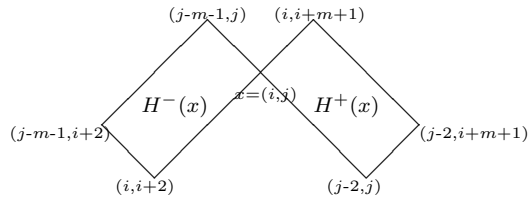
\[\begin{picture}(140,60)
\put(20,20){\line(1,-1){20}}
\put(-15,15){\tiny($j$-$m$-1,$i$+2)}
\put(70,30){\tiny$x$=($i$,$j$)}
\put(25,-5){\tiny($i$,$i$+2)}
\put(45,60){\tiny($j$-$m$-1,$j$)}
\put(35,25){\scriptsize$H^-(x)$}
\put(20,20){\line(1,1){40}}
\put(40,0){\line(1,1){60}}
\put(60,60){\line(1,-1){60}}
\put(120,0){\line(1,1){20}}
\put(100,60){\line(1,-1){40}}
\put(140,15){\tiny($j$-2,$i$+$m$+1)}
\put(85,60){\tiny($i$,$i$+$m$+1)}
\put(105,-5){\tiny($j$-2,$j$)}
\put(100,25){\scriptsize$H^+(x)$}
\end{picture}\]

If $\Delta=B_m$, $C_m$ or $D_{m+1}$, then $H^{\pm}(x)$ consists of points in the polygons (and their boundary) below.
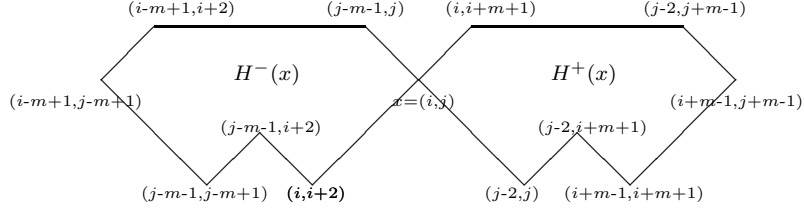
\[\begin{picture}(260,70)
\put(20,40){\line(1,1){20}}
\put(20,40){\line(1,-1){40}}
\put(60,0){\line(1,1){20}}
\put(40,60){\line(1,0){80}}
\put(80,20){\line(1,-1){20}}
\put(100,0){\line(1,1){60}}
\put(120,60){\line(1,-1){60}}
\put(130,30){\tiny$x$=($i$,$j$)}
\put(90,-5){\tiny($i$,$i$+2)}
\put(35,-5){\tiny($j$-$m$-1,$j$-$m$+1)}
\put(105,65){\tiny($j$-$m$-1,$j$)}
\put(30,65){\tiny($i$-$m$+1,$i$+2)}
\put(-15,30){\tiny($i$-$m$+1,$j$-$m$+1)}
\put(65,20){\tiny($j$-$m$-1,$i$+2)}
\put(70,40){\scriptsize$H^-(x)$}
\put(180,0){\line(1,1){20}}
\put(160,60){\line(1,0){80}}
\put(200,20){\line(1,-1){20}}
\put(220,0){\line(1,1){40}}
\put(240,60){\line(1,-1){20}}
\put(90,-5){\tiny($i$,$i$+2)}
\put(165,-5){\tiny($j$-2,$j$)}
\put(195,-5){\tiny($i$+$m$-1,$i$+$m$+1)}
\put(225,65){\tiny($j$-2,$j$+$m$-1)}
\put(150,65){\tiny($i$,$i$+$m$+1)}
\put(235,30){\tiny($i$+$m$-1,$j$+$m$-1)}
\put(185,20){\tiny($j$-2,$i$+$m$+1)}
\put(190,40){\scriptsize$H^+(x)$}
\end{picture}\]

We need more explanation for the case $\Delta=D_{m+1}$. If $\sigma x=x$, then $H^-(x)$ (resp. $H^+(x)$) contains both of $(p,p+m+1)_+$ and $(p,p+m+1)_-$ for any $p$ satisfying $i-m+1\le p\le j-m-1$ (resp. $i\le p\le j-2$). If $\sigma x\neq x$, then $H^{\pm}(x)$ contains $(\tau\sigma)^{\mp p}x$ ($0\le p<m$), but does not contain $\sigma(\tau\sigma)^{\mp p}x$ ($0\le p<m$).

Define an automorphism $\omega$ on $\zzz \Delta$ by $\omega(i,j):=(j-m-2,i+1)$ for $\Delta=A_m$, and $\omega:=\sigma(\tau\sigma)^m$ for other $\Delta$. It is easily checked that $\tau$, $\omega$ and $H^{\pm}$ mutually commute. Put $\tau_i:=\tau\omega^{i-1}$. We call a subset $S$ of $(\zzz\Delta)_0$ {\it maximal $n$-orthogonal} if
\[(\zzz\Delta)_0\backslash S=\bigcup_{x\in S,\ 0<i\le n}H^-(\tau_ix).\]

In \XDDC, we will show that these combinatorial datum $\omega$ and $H^-(x)$ (resp. $H^+(x)$) on $\zzz\Delta$ correspond to homological datum $\Omega$ and $\{Y\in\ind{}_\Lambda\underline{\dn{M}}\ |\ \underline{\hom}_\Lambda(Y,X)\neq0$ (resp. $\underline{\hom}_\Lambda(X,Y)\neq0$)$\}$ on ${}_\Lambda\underline{\dn{M}}$.

\vskip.5em{\bf\XDBA\ Proposition }{\it
(1) $x\in H^-(y)$ if and only if $y\in H^+(x)$.

(2) $H^-(\tau x)=H^+(\omega x)$.

(3) Any maximal $n$-orthogonal subset $S$ satisfies $\tau_{n+1}S=S$.

(4) $S$ is maximal $n$-orthogonal if and only if $(\zzz\Delta)_0\backslash S=\bigcup_{x\in S,\ 0<i\le n}H^+(\tau_i^{-1}x)$.

(5) Any $\tau_{n+1}$-invariant subset $S$ of $(\zzz\Delta)_0$ satisfying $S\cap(\bigcup_{x\in S,\ 0<i\le n}H^-(\tau_ix))=\emptyset$ is contained in a maximal $n$-orthogonal subset.}

\vskip.5em{\sc Proof }
(1)(2) Easy.

(3) We will show that $\tau_{n+1}x\in H^-(\tau_iy)$ if and only if $y\in H^-(\tau_{n+1-i}x)$. Since $H^-(\tau_iy)=H^+(\omega^iy)$ holds by (2), $\tau_{n+1}x\in H^-(\tau_iy)$ if and only if $\omega^iy\in H^-(\tau_{n+1}x)$ by (1). Since $\omega^{-i}H^-(\tau_{n+1}x)=H^-(\tau_{n+1-i}x)$, the assertion follows.

If $x\in S$ and $\tau_{n+1}x\notin S$, then there exist $y\in S$ and $i$ ($0<i\le n$) such that $\tau_{n+1}x\in H^-(\tau_iy)$. Then $y\in H^-(\tau_{n+1-i}x)$ holds, a contradiction. If $x\notin S$ and $\tau_{n+1}x\in S$, then there exist $y\in S$ and $i$ ($0<i\le n$) such that $x\in H^-(\tau_iy)$. Then $\tau_{n+1}y\in S$ and $\tau_{n+1}x\in H^-(\tau_i\tau_{n+1}y)$ hold, a contradiction.

(4) Immediate from $H^-(\tau_ix)=H^+(\omega^ix)=H^+(\tau_{n+1-i}^{-1}\tau_{n+1}x)$ and (3).

(5) Assume $S$ is not maximal $n$-orthogonal. Take $y\in(\zzz\Delta)_0$ which is not contained in $S\cup(\bigcup_{x\in S,\ 0<i\le n}H^-(\tau_ix))$. Then $S^\prime:=S\cup\{\tau_{n+1}^iy\ |\ i\in\zzz\}$ satisfies $S^\prime\cap(\bigcup_{x\in S^\prime,\ 0<i\le n}H^-(\tau_ix))=\emptyset$ by the proof of (3). Repeating similar argument, we obtain the assertion.\rule{5pt}{10pt}

\vskip.5em{\bf\XDBB\ }We call a subset $S$ of $\coprod_k(\zzz\Delta_k)_0$ {\it maximal $n$-orthogonal} if $S\cup(\zzz\Delta_k)_0$ is a maximal $n$-orthogonal subset of $(\zzz\Delta_k)_0$ for any $k$. We will prove the theorems below in \XDD.

\vskip.5em{\bf Theorem }{\it
Let $\dn{A}({}_\Lambda\underline{\dn{M}})=\coprod_k(\zzz\Delta_k/G_k)$ and $p:(\zzz\Delta_k)_0\to(\zzz\Delta_k/G_k)_0=\ind{}_\Lambda\underline{\dn{M}}$ the surjection. Let $\cc$ be a subcategory of ${}_\Lambda\dn{M}$ such that $\Lambda\in\cc$. Then $\cc$ is maximal $n$-orthogonal if and only if the subset $p^{-1}(\ind\underline{\cc})$ of $\coprod_k(\zzz\Delta_k)_0$ is maximal $n$-orthogonal.}

\vskip.5em{\bf\XDBC\ Theorem }{\it
The number of maximal $1$-orthogonal subsets of $(\zzz\Delta)_0$ is $\frac{1}{m+2}{2m+2\choose m+1}$ if $\Delta=A_m$, ${2m\choose m}$ if $\Delta=B_m$ or $C_m$, and $\frac{3m+1}{m+1}{2m\choose m}$ if $\Delta=D_{m+1}$.}

\vskip.5em{\bf\XDBD\ Remark }
We notice that \XDBB\ looks like Riedtmann and Wiedemann's classification of configurations of projective modules for selfinjective algebras and one-dimensional Gorenstein orders $\Lambda$ [R1][W1,2]. Their method heavily depends on the assumption that $\Lambda$ is finite dimensional over an algebraically closed field. To treat arbitrary $\Lambda$ in \XDD, we will need results in [I1,2,3] which are obtained by functorial methods.

\vskip.5em{\bf\XDBE\ Example }
We will indicate all maximal $1$-orthogonal subsets $S$ of $(\zzz\Delta)_0$ up to the automorphism group $\langle\tau,\sigma\rangle$ of $\zzz\Delta$. We will encircle vertices contained in $S$. We enumerate the number of maximal $1$-orthogonal subsets of $(\zzz\Delta)_0$.
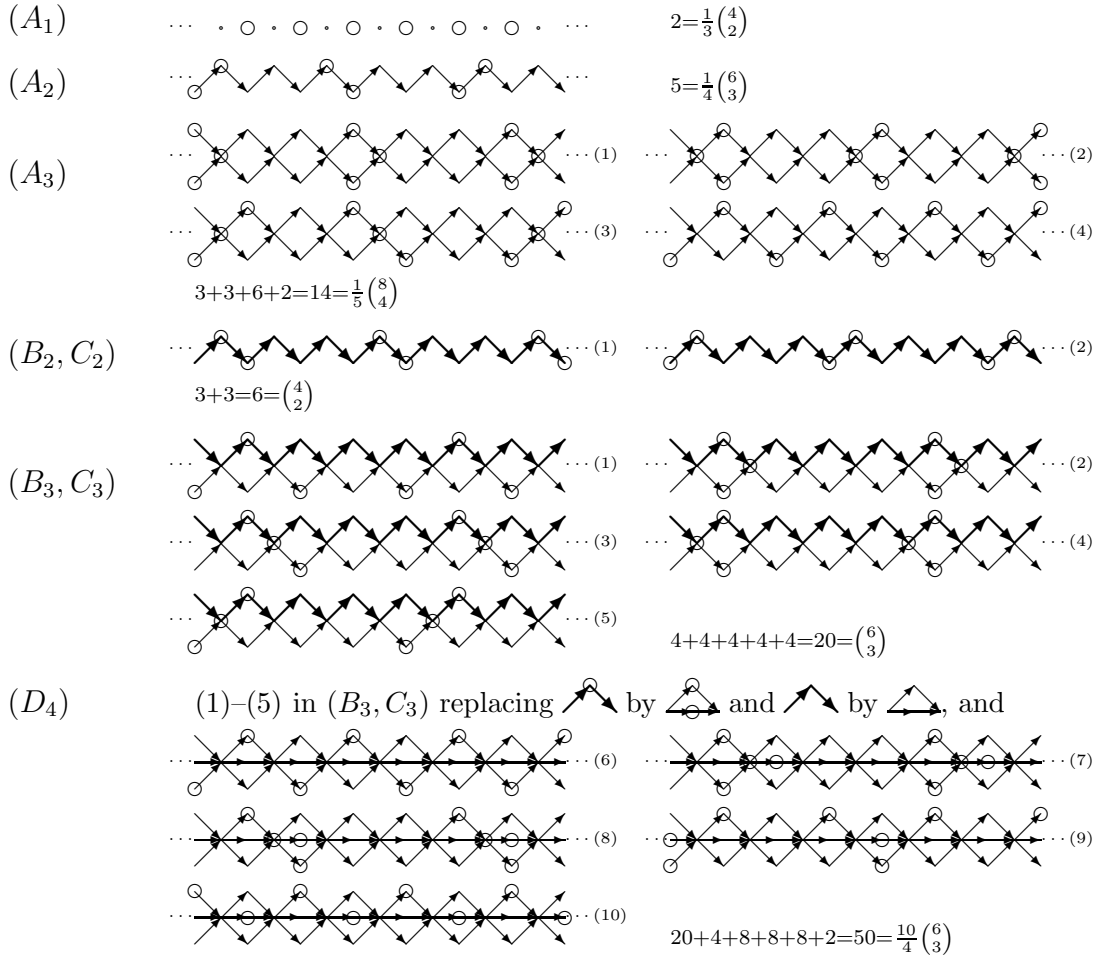
\[\begin{array}{ll}
(A_1)&\begin{picture}(180,5)
\put(10,0){\circle{1}}\put(20,0){\circle{5}}\put(30,0){\circle{1}}\put(40,0){\circle{5}}\put(50,0){\circle{1}}\put(60,0){\circle{5}}\put(70,0){\circle{1}}\put(80,0){\circle{5}}\put(90,0){\circle{1}}\put(100,0){\circle{5}}\put(110,0){\circle{1}}\put(120,0){\circle{5}}\put(130,0){\circle{1}}
\put(140,-1){\tiny$\cdots$}\put(-10,-1){\tiny$\cdots$}
\end{picture}{\scriptstyle 2=\frac{1}{3}{4\choose2}}\\
(A_2)&\begin{picture}(180,20)
\put(0,0){\vector(1,1){10}}\put(10,10){\vector(1,-1){10}}\put(20,0){\vector(1,1){10}}\put(30,10){\vector(1,-1){10}}\put(40,0){\vector(1,1){10}}\put(50,10){\vector(1,-1){10}}\put(60,0){\vector(1,1){10}}\put(70,10){\vector(1,-1){10}}\put(80,0){\vector(1,1){10}}\put(90,10){\vector(1,-1){10}}\put(100,0){\vector(1,1){10}}\put(110,10){\vector(1,-1){10}}\put(120,0){\vector(1,1){10}}\put(130,10){\vector(1,-1){10}}
\put(0,0){\circle{5}}\put(10,10){\circle{5}}\put(50,10){\circle{5}}\put(60,0){\circle{5}}\put(100,0){\circle{5}}\put(110,10){\circle{5}}
\put(140,4){\tiny$\cdots$}\put(-10,4){\tiny$\cdots$}
\end{picture}{\scriptstyle 5=\frac{1}{4}{6\choose3}}\\
(A_3)&\begin{picture}(180,30)
\put(0,0){\vector(1,1){10}}\put(10,10){\vector(1,-1){10}}\put(20,0){\vector(1,1){10}}\put(30,10){\vector(1,-1){10}}\put(40,0){\vector(1,1){10}}\put(50,10){\vector(1,-1){10}}\put(60,0){\vector(1,1){10}}\put(70,10){\vector(1,-1){10}}\put(80,0){\vector(1,1){10}}\put(90,10){\vector(1,-1){10}}\put(100,0){\vector(1,1){10}}\put(110,10){\vector(1,-1){10}}\put(120,0){\vector(1,1){10}}\put(130,10){\vector(1,-1){10}}
\put(0,20){\vector(1,-1){10}}\put(10,10){\vector(1,1){10}}\put(20,20){\vector(1,-1){10}}\put(30,10){\vector(1,1){10}}\put(40,20){\vector(1,-1){10}}\put(50,10){\vector(1,1){10}}\put(60,20){\vector(1,-1){10}}\put(70,10){\vector(1,1){10}}\put(80,20){\vector(1,-1){10}}\put(90,10){\vector(1,1){10}}\put(100,20){\vector(1,-1){10}}\put(110,10){\vector(1,1){10}}\put(120,20){\vector(1,-1){10}}\put(130,10){\vector(1,1){10}}
\put(0,0){\circle{5}}\put(10,10){\circle{5}}\put(0,20){\circle{5}}\put(60,20){\circle{5}}\put(70,10){\circle{5}}\put(60,0){\circle{5}}\put(120,0){\circle{5}}\put(130,10){\circle{5}}\put(120,20){\circle{5}}
\put(140,9){\tiny$\cdots(1)$}\put(-10,9){\tiny$\cdots$}
\end{picture}
\begin{picture}(180,30)
\put(0,0){\vector(1,1){10}}\put(10,10){\vector(1,-1){10}}\put(20,0){\vector(1,1){10}}\put(30,10){\vector(1,-1){10}}\put(40,0){\vector(1,1){10}}\put(50,10){\vector(1,-1){10}}\put(60,0){\vector(1,1){10}}\put(70,10){\vector(1,-1){10}}\put(80,0){\vector(1,1){10}}\put(90,10){\vector(1,-1){10}}\put(100,0){\vector(1,1){10}}\put(110,10){\vector(1,-1){10}}\put(120,0){\vector(1,1){10}}\put(130,10){\vector(1,-1){10}}
\put(0,20){\vector(1,-1){10}}\put(10,10){\vector(1,1){10}}\put(20,20){\vector(1,-1){10}}\put(30,10){\vector(1,1){10}}\put(40,20){\vector(1,-1){10}}\put(50,10){\vector(1,1){10}}\put(60,20){\vector(1,-1){10}}\put(70,10){\vector(1,1){10}}\put(80,20){\vector(1,-1){10}}\put(90,10){\vector(1,1){10}}\put(100,20){\vector(1,-1){10}}\put(110,10){\vector(1,1){10}}\put(120,20){\vector(1,-1){10}}\put(130,10){\vector(1,1){10}}
\put(20,0){\circle{5}}\put(10,10){\circle{5}}\put(20,20){\circle{5}}\put(80,20){\circle{5}}\put(70,10){\circle{5}}\put(80,0){\circle{5}}\put(140,0){\circle{5}}\put(130,10){\circle{5}}\put(140,20){\circle{5}}
\put(140,9){\tiny$\cdots(2)$}\put(-10,9){\tiny$\cdots$}
\end{picture}\\
&\begin{picture}(180,25)
\put(0,0){\vector(1,1){10}}\put(10,10){\vector(1,-1){10}}\put(20,0){\vector(1,1){10}}\put(30,10){\vector(1,-1){10}}\put(40,0){\vector(1,1){10}}\put(50,10){\vector(1,-1){10}}\put(60,0){\vector(1,1){10}}\put(70,10){\vector(1,-1){10}}\put(80,0){\vector(1,1){10}}\put(90,10){\vector(1,-1){10}}\put(100,0){\vector(1,1){10}}\put(110,10){\vector(1,-1){10}}\put(120,0){\vector(1,1){10}}\put(130,10){\vector(1,-1){10}}
\put(0,20){\vector(1,-1){10}}\put(10,10){\vector(1,1){10}}\put(20,20){\vector(1,-1){10}}\put(30,10){\vector(1,1){10}}\put(40,20){\vector(1,-1){10}}\put(50,10){\vector(1,1){10}}\put(60,20){\vector(1,-1){10}}\put(70,10){\vector(1,1){10}}\put(80,20){\vector(1,-1){10}}\put(90,10){\vector(1,1){10}}\put(100,20){\vector(1,-1){10}}\put(110,10){\vector(1,1){10}}\put(120,20){\vector(1,-1){10}}\put(130,10){\vector(1,1){10}}
\put(0,0){\circle{5}}\put(10,10){\circle{5}}\put(20,20){\circle{5}}\put(60,20){\circle{5}}\put(70,10){\circle{5}}\put(80,0){\circle{5}}\put(120,0){\circle{5}}\put(130,10){\circle{5}}\put(140,20){\circle{5}}
\put(140,9){\tiny$\cdots(3)$}\put(-10,9){\tiny$\cdots$}
\end{picture}
\begin{picture}(180,25)
\put(0,0){\vector(1,1){10}}\put(10,10){\vector(1,-1){10}}\put(20,0){\vector(1,1){10}}\put(30,10){\vector(1,-1){10}}\put(40,0){\vector(1,1){10}}\put(50,10){\vector(1,-1){10}}\put(60,0){\vector(1,1){10}}\put(70,10){\vector(1,-1){10}}\put(80,0){\vector(1,1){10}}\put(90,10){\vector(1,-1){10}}\put(100,0){\vector(1,1){10}}\put(110,10){\vector(1,-1){10}}\put(120,0){\vector(1,1){10}}\put(130,10){\vector(1,-1){10}}
\put(0,20){\vector(1,-1){10}}\put(10,10){\vector(1,1){10}}\put(20,20){\vector(1,-1){10}}\put(30,10){\vector(1,1){10}}\put(40,20){\vector(1,-1){10}}\put(50,10){\vector(1,1){10}}\put(60,20){\vector(1,-1){10}}\put(70,10){\vector(1,1){10}}\put(80,20){\vector(1,-1){10}}\put(90,10){\vector(1,1){10}}\put(100,20){\vector(1,-1){10}}\put(110,10){\vector(1,1){10}}\put(120,20){\vector(1,-1){10}}\put(130,10){\vector(1,1){10}}
\put(0,0){\circle{5}}\put(20,20){\circle{5}}\put(40,0){\circle{5}}\put(60,20){\circle{5}}\put(80,0){\circle{5}}\put(100,20){\circle{5}}\put(120,0){\circle{5}}\put(140,20){\circle{5}}
\put(140,9){\tiny$\cdots(4)$}\put(-10,9){\tiny$\cdots$}
\end{picture}\\
&{\scriptstyle 3+3+6+2=14=\frac{1}{5}{8\choose4}}\\
(B_2,C_2)&\begin{picture}(180,20)
\thicklines
\put(0,0){\vector(1,1){10}}\put(10,10){\vector(1,-1){10}}\put(20,0){\vector(1,1){10}}\put(30,10){\vector(1,-1){10}}\put(40,0){\vector(1,1){10}}\put(50,10){\vector(1,-1){10}}\put(60,0){\vector(1,1){10}}\put(70,10){\vector(1,-1){10}}\put(80,0){\vector(1,1){10}}\put(90,10){\vector(1,-1){10}}\put(100,0){\vector(1,1){10}}\put(110,10){\vector(1,-1){10}}\put(120,0){\vector(1,1){10}}\put(130,10){\vector(1,-1){10}}\thinlines
\put(20,0){\circle{5}}\put(10,10){\circle{5}}\put(70,10){\circle{5}}\put(80,0){\circle{5}}\put(140,0){\circle{5}}\put(130,10){\circle{5}}
\put(140,4){\tiny$\cdots(1)$}\put(-10,4){\tiny$\cdots$}
\end{picture}
\begin{picture}(180,20)
\thicklines
\put(0,0){\vector(1,1){10}}\put(10,10){\vector(1,-1){10}}\put(20,0){\vector(1,1){10}}\put(30,10){\vector(1,-1){10}}\put(40,0){\vector(1,1){10}}\put(50,10){\vector(1,-1){10}}\put(60,0){\vector(1,1){10}}\put(70,10){\vector(1,-1){10}}\put(80,0){\vector(1,1){10}}\put(90,10){\vector(1,-1){10}}\put(100,0){\vector(1,1){10}}\put(110,10){\vector(1,-1){10}}\put(120,0){\vector(1,1){10}}\put(130,10){\vector(1,-1){10}}\thinlines
\put(0,0){\circle{5}}\put(10,10){\circle{5}}\put(70,10){\circle{5}}\put(60,0){\circle{5}}\put(120,0){\circle{5}}\put(130,10){\circle{5}}
\put(140,4){\tiny$\cdots(2)$}\put(-10,4){\tiny$\cdots$}
\end{picture}\\
&{\scriptstyle 3+3=6={4\choose2}}\\
(B_3,C_3)\ \ \ \ \ &\begin{picture}(180,30)
\put(0,0){\vector(1,1){10}}\put(10,10){\vector(1,-1){10}}\put(20,0){\vector(1,1){10}}\put(30,10){\vector(1,-1){10}}\put(40,0){\vector(1,1){10}}\put(50,10){\vector(1,-1){10}}\put(60,0){\vector(1,1){10}}\put(70,10){\vector(1,-1){10}}\put(80,0){\vector(1,1){10}}\put(90,10){\vector(1,-1){10}}\put(100,0){\vector(1,1){10}}\put(110,10){\vector(1,-1){10}}\put(120,0){\vector(1,1){10}}\put(130,10){\vector(1,-1){10}}
\thicklines\put(0,20){\vector(1,-1){10}}\put(10,10){\vector(1,1){10}}\put(20,20){\vector(1,-1){10}}\put(30,10){\vector(1,1){10}}\put(40,20){\vector(1,-1){10}}\put(50,10){\vector(1,1){10}}\put(60,20){\vector(1,-1){10}}\put(70,10){\vector(1,1){10}}\put(80,20){\vector(1,-1){10}}\put(90,10){\vector(1,1){10}}\put(100,20){\vector(1,-1){10}}\put(110,10){\vector(1,1){10}}\put(120,20){\vector(1,-1){10}}\put(130,10){\vector(1,1){10}}\thinlines
\put(0,0){\circle{5}}\put(20,20){\circle{5}}\put(40,0){\circle{5}}\put(100,20){\circle{5}}\put(80,0){\circle{5}}\put(120,0){\circle{5}}
\put(140,9){\tiny$\cdots(1)$}\put(-10,9){\tiny$\cdots$}
\end{picture}
\begin{picture}(180,30)
\put(0,0){\vector(1,1){10}}\put(10,10){\vector(1,-1){10}}\put(20,0){\vector(1,1){10}}\put(30,10){\vector(1,-1){10}}\put(40,0){\vector(1,1){10}}\put(50,10){\vector(1,-1){10}}\put(60,0){\vector(1,1){10}}\put(70,10){\vector(1,-1){10}}\put(80,0){\vector(1,1){10}}\put(90,10){\vector(1,-1){10}}\put(100,0){\vector(1,1){10}}\put(110,10){\vector(1,-1){10}}\put(120,0){\vector(1,1){10}}\put(130,10){\vector(1,-1){10}}
\thicklines\put(0,20){\vector(1,-1){10}}\put(10,10){\vector(1,1){10}}\put(20,20){\vector(1,-1){10}}\put(30,10){\vector(1,1){10}}\put(40,20){\vector(1,-1){10}}\put(50,10){\vector(1,1){10}}\put(60,20){\vector(1,-1){10}}\put(70,10){\vector(1,1){10}}\put(80,20){\vector(1,-1){10}}\put(90,10){\vector(1,1){10}}\put(100,20){\vector(1,-1){10}}\put(110,10){\vector(1,1){10}}\put(120,20){\vector(1,-1){10}}\put(130,10){\vector(1,1){10}}\thinlines
\put(20,0){\circle{5}}\put(20,20){\circle{5}}\put(30,10){\circle{5}}\put(100,20){\circle{5}}\put(100,0){\circle{5}}\put(110,10){\circle{5}}
\put(140,9){\tiny$\cdots(2)$}\put(-10,9){\tiny$\cdots$}
\end{picture}\\
&\begin{picture}(180,25)
\put(0,0){\vector(1,1){10}}\put(10,10){\vector(1,-1){10}}\put(20,0){\vector(1,1){10}}\put(30,10){\vector(1,-1){10}}\put(40,0){\vector(1,1){10}}\put(50,10){\vector(1,-1){10}}\put(60,0){\vector(1,1){10}}\put(70,10){\vector(1,-1){10}}\put(80,0){\vector(1,1){10}}\put(90,10){\vector(1,-1){10}}\put(100,0){\vector(1,1){10}}\put(110,10){\vector(1,-1){10}}\put(120,0){\vector(1,1){10}}\put(130,10){\vector(1,-1){10}}
\thicklines\put(0,20){\vector(1,-1){10}}\put(10,10){\vector(1,1){10}}\put(20,20){\vector(1,-1){10}}\put(30,10){\vector(1,1){10}}\put(40,20){\vector(1,-1){10}}\put(50,10){\vector(1,1){10}}\put(60,20){\vector(1,-1){10}}\put(70,10){\vector(1,1){10}}\put(80,20){\vector(1,-1){10}}\put(90,10){\vector(1,1){10}}\put(100,20){\vector(1,-1){10}}\put(110,10){\vector(1,1){10}}\put(120,20){\vector(1,-1){10}}\put(130,10){\vector(1,1){10}}\thinlines
\put(40,0){\circle{5}}\put(20,20){\circle{5}}\put(30,10){\circle{5}}\put(100,20){\circle{5}}\put(120,0){\circle{5}}\put(110,10){\circle{5}}
\put(140,9){\tiny$\cdots(3)$}\put(-10,9){\tiny$\cdots$}
\end{picture}
\begin{picture}(180,25)
\put(0,0){\vector(1,1){10}}\put(10,10){\vector(1,-1){10}}\put(20,0){\vector(1,1){10}}\put(30,10){\vector(1,-1){10}}\put(40,0){\vector(1,1){10}}\put(50,10){\vector(1,-1){10}}\put(60,0){\vector(1,1){10}}\put(70,10){\vector(1,-1){10}}\put(80,0){\vector(1,1){10}}\put(90,10){\vector(1,-1){10}}\put(100,0){\vector(1,1){10}}\put(110,10){\vector(1,-1){10}}\put(120,0){\vector(1,1){10}}\put(130,10){\vector(1,-1){10}}
\thicklines\put(0,20){\vector(1,-1){10}}\put(10,10){\vector(1,1){10}}\put(20,20){\vector(1,-1){10}}\put(30,10){\vector(1,1){10}}\put(40,20){\vector(1,-1){10}}\put(50,10){\vector(1,1){10}}\put(60,20){\vector(1,-1){10}}\put(70,10){\vector(1,1){10}}\put(80,20){\vector(1,-1){10}}\put(90,10){\vector(1,1){10}}\put(100,20){\vector(1,-1){10}}\put(110,10){\vector(1,1){10}}\put(120,20){\vector(1,-1){10}}\put(130,10){\vector(1,1){10}}\thinlines
\put(20,0){\circle{5}}\put(20,20){\circle{5}}\put(10,10){\circle{5}}\put(100,20){\circle{5}}\put(100,0){\circle{5}}\put(90,10){\circle{5}}
\put(140,9){\tiny$\cdots(4)$}\put(-10,9){\tiny$\cdots$}
\end{picture}\\
&\begin{picture}(180,25)
\put(0,0){\vector(1,1){10}}\put(10,10){\vector(1,-1){10}}\put(20,0){\vector(1,1){10}}\put(30,10){\vector(1,-1){10}}\put(40,0){\vector(1,1){10}}\put(50,10){\vector(1,-1){10}}\put(60,0){\vector(1,1){10}}\put(70,10){\vector(1,-1){10}}\put(80,0){\vector(1,1){10}}\put(90,10){\vector(1,-1){10}}\put(100,0){\vector(1,1){10}}\put(110,10){\vector(1,-1){10}}\put(120,0){\vector(1,1){10}}\put(130,10){\vector(1,-1){10}}
\thicklines\put(0,20){\vector(1,-1){10}}\put(10,10){\vector(1,1){10}}\put(20,20){\vector(1,-1){10}}\put(30,10){\vector(1,1){10}}\put(40,20){\vector(1,-1){10}}\put(50,10){\vector(1,1){10}}\put(60,20){\vector(1,-1){10}}\put(70,10){\vector(1,1){10}}\put(80,20){\vector(1,-1){10}}\put(90,10){\vector(1,1){10}}\put(100,20){\vector(1,-1){10}}\put(110,10){\vector(1,1){10}}\put(120,20){\vector(1,-1){10}}\put(130,10){\vector(1,1){10}}\thinlines
\put(0,0){\circle{5}}\put(20,20){\circle{5}}\put(10,10){\circle{5}}\put(100,20){\circle{5}}\put(80,0){\circle{5}}\put(90,10){\circle{5}}
\put(140,9){\tiny$\cdots(5)$}\put(-10,9){\tiny$\cdots$}
\end{picture}{\scriptstyle 4+4+4+4+4=20={6\choose3}}\\
(D_4)&\mbox{\small (1)--(5) in $(B_3,C_3)$ replacing
\begin{picture}(20,20)\thicklines\put(0,0){\vector(1,1){10}}\put(10,10){\vector(1,-1){10}}\thinlines\put(10,10){\circle{5}}\end{picture} by
\begin{picture}(20,20)\put(0,0){\vector(1,1){10}}\put(10,10){\vector(1,-1){10}}\put(0,0){\vector(1,0){10}}\put(10,0){\vector(1,0){10}}\put(10,10){\circle{5}}\put(10,0){\circle{5}}\end{picture} and 
\begin{picture}(20,20)\thicklines\put(0,0){\vector(1,1){10}}\put(10,10){\vector(1,-1){10}}\thinlines\end{picture} by
\begin{picture}(20,20)\put(0,0){\vector(1,1){10}}\put(10,10){\vector(1,-1){10}}\put(0,0){\vector(1,0){10}}\put(10,0){\vector(1,0){10}}\end{picture}, and}\\
&\begin{picture}(180,25)
\put(0,0){\vector(1,1){10}}\put(10,10){\vector(1,-1){10}}\put(20,0){\vector(1,1){10}}\put(30,10){\vector(1,-1){10}}\put(40,0){\vector(1,1){10}}\put(50,10){\vector(1,-1){10}}\put(60,0){\vector(1,1){10}}\put(70,10){\vector(1,-1){10}}\put(80,0){\vector(1,1){10}}\put(90,10){\vector(1,-1){10}}\put(100,0){\vector(1,1){10}}\put(110,10){\vector(1,-1){10}}\put(120,0){\vector(1,1){10}}\put(130,10){\vector(1,-1){10}}
\put(0,10){\vector(1,0){10}}\put(10,10){\vector(1,0){10}}\put(20,10){\vector(1,0){10}}\put(30,10){\vector(1,0){10}}\put(40,10){\vector(1,0){10}}\put(50,10){\vector(1,0){10}}\put(60,10){\vector(1,0){10}}\put(70,10){\vector(1,0){10}}\put(80,10){\vector(1,0){10}}\put(90,10){\vector(1,0){10}}\put(100,10){\vector(1,0){10}}\put(110,10){\vector(1,0){10}}\put(120,10){\vector(1,0){10}}\put(130,10){\vector(1,0){10}}
\put(0,20){\vector(1,-1){10}}\put(10,10){\vector(1,1){10}}\put(20,20){\vector(1,-1){10}}\put(30,10){\vector(1,1){10}}\put(40,20){\vector(1,-1){10}}\put(50,10){\vector(1,1){10}}\put(60,20){\vector(1,-1){10}}\put(70,10){\vector(1,1){10}}\put(80,20){\vector(1,-1){10}}\put(90,10){\vector(1,1){10}}\put(100,20){\vector(1,-1){10}}\put(110,10){\vector(1,1){10}}\put(120,20){\vector(1,-1){10}}\put(130,10){\vector(1,1){10}}
\put(0,0){\circle{5}}\put(20,20){\circle{5}}\put(40,0){\circle{5}}\put(60,20){\circle{5}}\put(100,20){\circle{5}}\put(80,0){\circle{5}}\put(120,0){\circle{5}}\put(140,20){\circle{5}}
\put(140,9){\tiny$\cdots(6)$}\put(-10,9){\tiny$\cdots$}
\end{picture}
\begin{picture}(180,25)
\put(0,0){\vector(1,1){10}}\put(10,10){\vector(1,-1){10}}\put(20,0){\vector(1,1){10}}\put(30,10){\vector(1,-1){10}}\put(40,0){\vector(1,1){10}}\put(50,10){\vector(1,-1){10}}\put(60,0){\vector(1,1){10}}\put(70,10){\vector(1,-1){10}}\put(80,0){\vector(1,1){10}}\put(90,10){\vector(1,-1){10}}\put(100,0){\vector(1,1){10}}\put(110,10){\vector(1,-1){10}}\put(120,0){\vector(1,1){10}}\put(130,10){\vector(1,-1){10}}
\put(0,10){\vector(1,0){10}}\put(10,10){\vector(1,0){10}}\put(20,10){\vector(1,0){10}}\put(30,10){\vector(1,0){10}}\put(40,10){\vector(1,0){10}}\put(50,10){\vector(1,0){10}}\put(60,10){\vector(1,0){10}}\put(70,10){\vector(1,0){10}}\put(80,10){\vector(1,0){10}}\put(90,10){\vector(1,0){10}}\put(100,10){\vector(1,0){10}}\put(110,10){\vector(1,0){10}}\put(120,10){\vector(1,0){10}}\put(130,10){\vector(1,0){10}}
\put(0,20){\vector(1,-1){10}}\put(10,10){\vector(1,1){10}}\put(20,20){\vector(1,-1){10}}\put(30,10){\vector(1,1){10}}\put(40,20){\vector(1,-1){10}}\put(50,10){\vector(1,1){10}}\put(60,20){\vector(1,-1){10}}\put(70,10){\vector(1,1){10}}\put(80,20){\vector(1,-1){10}}\put(90,10){\vector(1,1){10}}\put(100,20){\vector(1,-1){10}}\put(110,10){\vector(1,1){10}}\put(120,20){\vector(1,-1){10}}\put(130,10){\vector(1,1){10}}
\put(20,0){\circle{5}}\put(20,20){\circle{5}}\put(30,10){\circle{5}}\put(40,10){\circle{5}}\put(100,20){\circle{5}}\put(100,0){\circle{5}}\put(110,10){\circle{5}}\put(120,10){\circle{5}}
\put(140,9){\tiny$\cdots(7)$}\put(-10,9){\tiny$\cdots$}
\end{picture}\\
&\begin{picture}(180,25)
\put(0,0){\vector(1,1){10}}\put(10,10){\vector(1,-1){10}}\put(20,0){\vector(1,1){10}}\put(30,10){\vector(1,-1){10}}\put(40,0){\vector(1,1){10}}\put(50,10){\vector(1,-1){10}}\put(60,0){\vector(1,1){10}}\put(70,10){\vector(1,-1){10}}\put(80,0){\vector(1,1){10}}\put(90,10){\vector(1,-1){10}}\put(100,0){\vector(1,1){10}}\put(110,10){\vector(1,-1){10}}\put(120,0){\vector(1,1){10}}\put(130,10){\vector(1,-1){10}}
\put(0,10){\vector(1,0){10}}\put(10,10){\vector(1,0){10}}\put(20,10){\vector(1,0){10}}\put(30,10){\vector(1,0){10}}\put(40,10){\vector(1,0){10}}\put(50,10){\vector(1,0){10}}\put(60,10){\vector(1,0){10}}\put(70,10){\vector(1,0){10}}\put(80,10){\vector(1,0){10}}\put(90,10){\vector(1,0){10}}\put(100,10){\vector(1,0){10}}\put(110,10){\vector(1,0){10}}\put(120,10){\vector(1,0){10}}\put(130,10){\vector(1,0){10}}
\put(0,20){\vector(1,-1){10}}\put(10,10){\vector(1,1){10}}\put(20,20){\vector(1,-1){10}}\put(30,10){\vector(1,1){10}}\put(40,20){\vector(1,-1){10}}\put(50,10){\vector(1,1){10}}\put(60,20){\vector(1,-1){10}}\put(70,10){\vector(1,1){10}}\put(80,20){\vector(1,-1){10}}\put(90,10){\vector(1,1){10}}\put(100,20){\vector(1,-1){10}}\put(110,10){\vector(1,1){10}}\put(120,20){\vector(1,-1){10}}\put(130,10){\vector(1,1){10}}
\put(20,20){\circle{5}}\put(30,10){\circle{5}}\put(40,10){\circle{5}}\put(40,0){\circle{5}}\put(100,20){\circle{5}}\put(110,10){\circle{5}}\put(120,10){\circle{5}}\put(120,0){\circle{5}}
\put(140,9){\tiny$\cdots(8)$}\put(-10,9){\tiny$\cdots$}
\end{picture}
\begin{picture}(180,25)
\put(0,0){\vector(1,1){10}}\put(10,10){\vector(1,-1){10}}\put(20,0){\vector(1,1){10}}\put(30,10){\vector(1,-1){10}}\put(40,0){\vector(1,1){10}}\put(50,10){\vector(1,-1){10}}\put(60,0){\vector(1,1){10}}\put(70,10){\vector(1,-1){10}}\put(80,0){\vector(1,1){10}}\put(90,10){\vector(1,-1){10}}\put(100,0){\vector(1,1){10}}\put(110,10){\vector(1,-1){10}}\put(120,0){\vector(1,1){10}}\put(130,10){\vector(1,-1){10}}
\put(0,10){\vector(1,0){10}}\put(10,10){\vector(1,0){10}}\put(20,10){\vector(1,0){10}}\put(30,10){\vector(1,0){10}}\put(40,10){\vector(1,0){10}}\put(50,10){\vector(1,0){10}}\put(60,10){\vector(1,0){10}}\put(70,10){\vector(1,0){10}}\put(80,10){\vector(1,0){10}}\put(90,10){\vector(1,0){10}}\put(100,10){\vector(1,0){10}}\put(110,10){\vector(1,0){10}}\put(120,10){\vector(1,0){10}}\put(130,10){\vector(1,0){10}}
\put(0,20){\vector(1,-1){10}}\put(10,10){\vector(1,1){10}}\put(20,20){\vector(1,-1){10}}\put(30,10){\vector(1,1){10}}\put(40,20){\vector(1,-1){10}}\put(50,10){\vector(1,1){10}}\put(60,20){\vector(1,-1){10}}\put(70,10){\vector(1,1){10}}\put(80,20){\vector(1,-1){10}}\put(90,10){\vector(1,1){10}}\put(100,20){\vector(1,-1){10}}\put(110,10){\vector(1,1){10}}\put(120,20){\vector(1,-1){10}}\put(130,10){\vector(1,1){10}}
\put(0,0){\circle{5}}\put(0,10){\circle{5}}\put(20,20){\circle{5}}\put(60,20){\circle{5}}\put(80,0){\circle{5}}\put(80,10){\circle{5}}\put(100,20){\circle{5}}\put(140,20){\circle{5}}
\put(140,9){\tiny$\cdots(9)$}\put(-10,9){\tiny$\cdots$}\end{picture}\\
&\begin{picture}(180,25)\put(0,0){\vector(1,1){10}}\put(10,10){\vector(1,-1){10}}\put(20,0){\vector(1,1){10}}\put(30,10){\vector(1,-1){10}}\put(40,0){\vector(1,1){10}}\put(50,10){\vector(1,-1){10}}\put(60,0){\vector(1,1){10}}\put(70,10){\vector(1,-1){10}}\put(80,0){\vector(1,1){10}}\put(90,10){\vector(1,-1){10}}\put(100,0){\vector(1,1){10}}\put(110,10){\vector(1,-1){10}}\put(120,0){\vector(1,1){10}}\put(130,10){\vector(1,-1){10}} \put(0,10){\vector(1,0){10}}\put(10,10){\vector(1,0){10}}\put(20,10){\vector(1,0){10}}\put(30,10){\vector(1,0){10}}\put(40,10){\vector(1,0){10}}\put(50,10){\vector(1,0){10}}\put(60,10){\vector(1,0){10}}\put(70,10){\vector(1,0){10}}\put(80,10){\vector(1,0){10}}\put(90,10){\vector(1,0){10}}\put(100,10){\vector(1,0){10}}\put(110,10){\vector(1,0){10}}\put(120,10){\vector(1,0){10}}\put(130,10){\vector(1,0){10}} \put(0,20){\vector(1,-1){10}}\put(10,10){\vector(1,1){10}}\put(20,20){\vector(1,-1){10}}\put(30,10){\vector(1,1){10}}\put(40,20){\vector(1,-1){10}}\put(50,10){\vector(1,1){10}}\put(60,20){\vector(1,-1){10}}\put(70,10){\vector(1,1){10}}\put(80,20){\vector(1,-1){10}}\put(90,10){\vector(1,1){10}}\put(100,20){\vector(1,-1){10}}\put(110,10){\vector(1,1){10}}\put(120,20){\vector(1,-1){10}}\put(130,10){\vector(1,1){10}} \put(0,20){\circle{5}}\put(20,10){\circle{5}}\put(40,20){\circle{5}}\put(60,10){\circle{5}}\put(80,20){\circle{5}}\put(100,10){\circle{5}}\put(120,20){\circle{5}}\put(140,10){\circle{5}} \put(140,9){\tiny$\cdots(10)$}\put(-10,9){\tiny$\cdots$}\end{picture}{\scriptstyle 20+4+8+8+8+2=50=\frac{10}{4}{6\choose3}}
\end{array}\]

\vskip.5em{\bf\XDC\ }
Put $l:=m+3$ if $\Delta=A_m$, and $l:=2m+2$ if $\Delta=B_m$, $C_m$ or $D_{m+1}$. Then $\tau_2^2=\tau^l$ holds. Since any maximal $1$-orthogonal subset $S$ of $(\zzz\Delta)_0$ is $\tau_2$-invariant by \XDBA(3), $S$ can be regarded as a subset of $(\zzz\Delta/\langle\tau^l\rangle)_0=((\zzz/l\zzz)\Delta)_0$. We will show that such subsets can be displayed very clearly.

\vskip.5em{\bf\XDCA\ Definition }
Let $R_{l}$ be a regular $l$-polygon and $D(R_{l})$ the set of all diagonals of $R_{l}$ except edges of $R_{l}$. Consider the following conditions on a subset $S$ of $D(R_{l})$.

\strut\kern1em(i) Two diagonals in $S$ do not cross except their endopoints.

\strut\kern1em(i)$^\prime$ Two diagonals in $S$ do not cross except their endopoints and the center of $R_{l}$.

\strut\kern1em(ii) $R_{l}$ is dissected into triangles by diagonals in $S$.

\strut\kern1em(iii) $S$ is symmetric with respect to the center of $R_{l}$.

We call a subset $S$ of $D(R_{l})$ {\it maximal $1$-orthogonal} if it satisfies the condition (i) and (ii) (if $\Delta=A_m$), (i),(ii) and (iii) (if $\Delta=B_m$ or $C_m$), (i)$^\prime$,(ii) and (iii) (if $\Delta=D_{m+1}$).
In \XDCE, we will show that maximal $1$-orthogonal subsets of $(\zzz\Delta)_0$ corresponds to maximal $1$-orthogonal subsets of $D(R_{l})$. The number of maximal $1$-orthogonal subsets of $D(R_{l})$ is well-known for the case $\Delta=A_m$ [St], and the number for the case $\Delta=D_{m+1}$ is slightly different from \XDBC.

\vskip.5em{\bf\XDCB\ Theorem }{\it
The number of maximal $1$-orthogonal subsets of $D(R_{l})$ is $\frac{1}{m+2}{2m+2\choose m+1}$ if $\Delta=A_m$, ${2m\choose m}$ if $\Delta=B_m$ or $C_m$, and ${2m+1\choose m}$ if $\Delta=D_{m+1}$.}

\vskip.5em{\sc Proof }
It is well-known that the number of ways to dissect $R_{l}$ into triangles by diagonals without their crossing is given by Catalan number $a_l:=\frac{1}{l-1}{2l-4\choose l-2}$ [St]. This immediately implies the equality for $A_m$. Let $\Delta=B_m$ or $C_m$. Then any maximal $1$-orthogonal subset contains only one main diagonal of $R_{l}$. Of course, any main diagonal dissects $R_{l}$ into two congruent $(m+2)$-polygons. Since the number of main diagonals of $R_l$ is $m+1$ and the number of desired partitions of a $(m+2)$-polygon is $a_{m+2}=\frac{1}{m+1}{2m\choose m}$, we obtain the equality for $B_m$ and $C_m$. We will show the assertion for $D_{m+1}$. It is well-known that $f(x):=\sum_{i\ge0}a_{i+2}x^i$ is given by $f(x)=\frac{1-\sqrt{1-4x}}{2x}$ [St]. Put $b_0:=\frac{1}{2}$, $b_1:=1$, and for $i\ge2$, let $b_m$ be the number of the subsets $S$ of $D(R_{2m})$ satisfying the conditions (i)$^\prime$,(ii) and (iii) in \XDCA. We will show the recursion formula $b_{m}=2\sum_{i=2}^{m+1}a_ib_{m-i+1}$, which implies that $g(x):=\sum_{i\ge0}b_ix^i$ satisfies $2xfg=g-\frac{1}{2}$. Then $g=\frac{1}{2(1-2xf)}=\frac{1}{2\sqrt{1-4x}}=\sum_{i\ge0}{2i-1\choose i}x^i$ shows $b_i={2i-1\choose i}$. We denote by $c$ the center of $R_{2m}$. For any subset $S$ satisfying (i)$^\prime$,(ii) and (iii) in \XDCA, there exists a unique $v\in\{c,2,3,\cdots,2m-1\}$ such that there is some triangle $T_v$ with the set $\{0,1,v\}$ of vertices. By (iii), there is another triangle $T_v^\prime$ with the set $\{m,m+1,v+m\}$ ($v\neq c$) or $\{m,m+1,c\}$ ($v=c$) of vertices. Let $b_{m,v}$ be the number of such $S$. Then $b_m=b_{m,c}+\sum_{v=2}^{2m-1}b_{m,v}$ holds. If $v\in\{2,\cdots,m\}$, then $R_{2m}$ is dissected into two symmetric $v$-polygons and a $2(m-v+1)$-polygon by $T_v$ and $T_v^\prime$, and it is easily checked that $b_{m,v}=a_vb_{m-v+1}$ holds. If $v\in\{m+1,\cdots,2m-1\}$, then $b_{m,v}=b_{m,2m+1-v}$ holds by symmetry. If $v=c$, then $R_{2m}$ is dissected into two symmetric $(m+1)$-polygons by $T_v$ and $T_v^\prime$, and it is easily checked that $b_{m,c}=a_{m+1}$ holds. Thus $b_m=a_{m+1}+2\sum_{v=2}^{m}a_vb_{m-v+1}=2\sum_{v=2}^{m+1}a_vb_{m-v+1}$ holds.\rule{5pt}{10pt}

\vskip.5em{\bf\XDCC\ Example }
We indicate all maximal $1$-orthogonal subsets of $D(R_{l})$ up to rotations of $R_{l}$. In \XDCE\ below, we will give a correspondence between maximal $1$-orthogonal subsets of $D(R_{l})$ and those of $(\zzz\Delta)_0$. The correspondence between the examples below and those in \XDBE\ is indicated by indices. We enumerate the number without identifying up to rotations.
\[\begin{array}{lcr}
(A_1)&\begin{picture}(30,20)\put(0,10){\line(1,-1){10}}\put(0,10){\line(-1,-1){10}}\put(0,-10){\line(1,1){10}}\put(0,-10){\line(-1,1){10}}\put(0,10){\line(0,-1){20}}\end{picture}&{\scriptstyle 2=\frac{1}{3}{4\choose2}}\\
(A_2)\ \ \ \ \ \ &\begin{picture}(30,20)
\put(0,10){\line(1,-1){10}}\put(0,10){\line(-1,-1){10}}\put(10,0){\line(-1,-2){5}}\put(-10,0){\line(1,-2){5}}\put(5,-10){\line(-1,0){10}}
\put(0,10){\line(1,-4){5}}\put(0,10){\line(-1,-4){5}}
\end{picture}&{\scriptstyle 5=\frac{1}{4}{6\choose3}}\\
(A_3)\ \ \ \ \ \ &
\begin{picture}(35,20)
\put(10,5){\line(0,-1){10}}\put(-10,5){\line(0,-1){10}}\put(0,10){\line(2,-1){10}}\put(0,10){\line(-2,-1){10}}\put(0,-10){\line(2,1){10}}\put(0,-10){\line(-2,1){10}}
\put(0,10){\line(0,-1){20}}\put(0,10){\line(-2,-3){10}}\put(0,-10){\line(2,3){10}}\put(12,-2){\tiny$(1)$}
\end{picture}
\begin{picture}(35,20)
\put(10,5){\line(0,-1){10}}\put(-10,5){\line(0,-1){10}}\put(0,10){\line(2,-1){10}}\put(0,10){\line(-2,-1){10}}\put(0,-10){\line(2,1){10}}\put(0,-10){\line(-2,1){10}}
\put(0,10){\line(0,-1){20}}\put(0,-10){\line(-2,3){10}}\put(0,10){\line(2,-3){10}}\put(12,-2){\tiny$(2)$}
\end{picture}
\begin{picture}(35,20)
\put(10,5){\line(0,-1){10}}\put(-10,5){\line(0,-1){10}}\put(0,10){\line(2,-1){10}}\put(0,10){\line(-2,-1){10}}\put(0,-10){\line(2,1){10}}\put(0,-10){\line(-2,1){10}}
\put(0,10){\line(0,-1){20}}\put(0,10){\line(-2,-3){10}}\put(0,10){\line(2,-3){10}}\put(12,-2){\tiny$(3)$}
\end{picture}
\begin{picture}(35,20)
\put(10,5){\line(0,-1){10}}\put(-10,5){\line(0,-1){10}}\put(0,10){\line(2,-1){10}}\put(0,10){\line(-2,-1){10}}\put(0,-10){\line(2,1){10}}\put(0,-10){\line(-2,1){10}}
\put(0,10){\line(2,-3){10}}\put(0,10){\line(-2,-3){10}}\put(-10,-5){\line(1,0){20}}\put(12,-2){\tiny$(4)$}
\end{picture}&{\scriptstyle 3+3+6+2=14=\frac{1}{5}{8\choose4}}\\
(B_2,C_2)\ \ \ \ \ \ &\begin{picture}(35,20)
\put(10,5){\line(0,-1){10}}\put(-10,5){\line(0,-1){10}}\put(0,10){\line(2,-1){10}}\put(0,10){\line(-2,-1){10}}\put(0,-10){\line(2,1){10}}\put(0,-10){\line(-2,1){10}}
\put(0,10){\line(0,-1){20}}\put(0,-10){\line(-2,3){10}}\put(0,10){\line(2,-3){10}}\put(12,-2){\tiny$(1)$}
\end{picture}
\begin{picture}(35,20)
\put(10,5){\line(0,-1){10}}\put(-10,5){\line(0,-1){10}}\put(0,10){\line(2,-1){10}}\put(0,10){\line(-2,-1){10}}\put(0,-10){\line(2,1){10}}\put(0,-10){\line(-2,1){10}}
\put(0,10){\line(0,-1){20}}\put(0,10){\line(-2,-3){10}}\put(0,-10){\line(2,3){10}}\put(12,-2){\tiny$(2)$}
\end{picture}&{\scriptstyle 3+3=6={4\choose2}}\\
(B_3,C_3)\ \ \ \ \ \ \ \ &\begin{picture}(45,30)
\put(15,0){\line(-1,2){5}}\put(15,0){\line(-1,-2){5}}\put(0,15){\line(2,-1){10}}\put(0,15){\line(-2,-1){10}}\put(-15,0){\line(1,2){5}}\put(-15,0){\line(1,-2){5}}\put(0,-15){\line(2,1){10}}\put(0,-15){\line(-2,1){10}}
\put(0,15){\line(0,-1){30}}\put(0,15){\line(-1,-1){15}}\put(0,15){\line(1,-1){15}}\put(0,-15){\line(1,1){15}}\put(0,-15){\line(-1,1){15}}\put(17,-2){\tiny$(1)$}
\end{picture}
\begin{picture}(45,30)
\put(15,0){\line(-1,2){5}}\put(15,0){\line(-1,-2){5}}\put(0,15){\line(2,-1){10}}\put(0,15){\line(-2,-1){10}}\put(-15,0){\line(1,2){5}}\put(-15,0){\line(1,-2){5}}\put(0,-15){\line(2,1){10}}\put(0,-15){\line(-2,1){10}}
\put(0,15){\line(0,-1){30}}\put(0,15){\line(2,-5){10}}\put(-10,-10){\line(0,1){20}}\put(0,-15){\line(-2,5){10}}\put(10,-10){\line(0,1){20}}\put(17,-2){\tiny$(2)$}
\end{picture}
\begin{picture}(45,30)
\put(15,0){\line(-1,2){5}}\put(15,0){\line(-1,-2){5}}\put(0,15){\line(2,-1){10}}\put(0,15){\line(-2,-1){10}}\put(-15,0){\line(1,2){5}}\put(-15,0){\line(1,-2){5}}\put(0,-15){\line(2,1){10}}\put(0,-15){\line(-2,1){10}}
\put(0,15){\line(0,-1){30}}\put(0,15){\line(2,-5){10}}\put(0,15){\line(1,-1){15}}\put(0,-15){\line(-2,5){10}}\put(0,-15){\line(-1,1){15}}\put(17,-2){\tiny$(3)$}
\end{picture}
\begin{picture}(45,30)
\put(15,0){\line(-1,2){5}}\put(15,0){\line(-1,-2){5}}\put(0,15){\line(2,-1){10}}\put(0,15){\line(-2,-1){10}}\put(-15,0){\line(1,2){5}}\put(-15,0){\line(1,-2){5}}\put(0,-15){\line(2,1){10}}\put(0,-15){\line(-2,1){10}}
\put(0,15){\line(0,-1){30}}\put(0,15){\line(-2,-5){10}}\put(-10,-10){\line(0,1){20}}\put(0,-15){\line(2,5){10}}\put(10,-10){\line(0,1){20}}\put(17,-2){\tiny$(4)$}
\end{picture}
\begin{picture}(15,30)
\put(15,0){\line(-1,2){5}}\put(15,0){\line(-1,-2){5}}\put(0,15){\line(2,-1){10}}\put(0,15){\line(-2,-1){10}}\put(-15,0){\line(1,2){5}}\put(-15,0){\line(1,-2){5}}\put(0,-15){\line(2,1){10}}\put(0,-15){\line(-2,1){10}}
\put(0,15){\line(0,-1){30}}\put(0,15){\line(-1,-1){15}}\put(0,15){\line(-2,-5){10}}\put(0,-15){\line(1,1){15}}\put(0,-15){\line(2,5){10}}\put(17,-2){\tiny$(5)$}
\end{picture}&{\scriptstyle 4+4+4+4+4=20={6\choose3}}\\
(D_3)\ \ \ \ \ \ &
\mbox{\small (1) and (2) in $(B_2,C_2)$ and }\ \ \ \ \ 
\begin{picture}(35,30)
\put(10,-5){\line(-2,1){20}}
\put(10,5){\line(0,-1){10}}\put(-10,5){\line(0,-1){10}}\put(0,10){\line(2,-1){10}}\put(0,10){\line(-2,-1){10}}\put(0,-10){\line(2,1){10}}\put(0,-10){\line(-2,1){10}}
\put(0,10){\line(0,-1){20}}\put(0,-10){\line(-2,3){10}}\put(0,10){\line(2,-3){10}}\put(12,-2){\tiny$(3)$}
\end{picture}
\begin{picture}(35,30)
\put(10,-5){\line(-2,1){20}}\put(10,5){\line(-2,-1){20}}
\put(10,5){\line(0,-1){10}}\put(-10,5){\line(0,-1){10}}\put(0,10){\line(2,-1){10}}\put(0,10){\line(-2,-1){10}}\put(0,-10){\line(2,1){10}}\put(0,-10){\line(-2,1){10}}
\put(0,10){\line(0,-1){20}}\put(12,-2){\tiny$(4)$}
\end{picture}&{\scriptstyle 6+3+1=10={5\choose2}}\\
(D_4)\ \ \ \ \ \ \ \ \begin{picture}(5,20)\end{picture}&\mbox{\small (1)--(5) in $(B_3,C_3)$ and }&\\
&\begin{picture}(45,20)
\put(15,0){\line(-1,2){5}}\put(15,0){\line(-1,-2){5}}\put(0,15){\line(2,-1){10}}\put(0,15){\line(-2,-1){10}}\put(-15,0){\line(1,2){5}}\put(-15,0){\line(1,-2){5}}\put(0,-15){\line(2,1){10}}\put(0,-15){\line(-2,1){10}}
\put(0,15){\line(0,-1){30}}\put(0,15){\line(-1,-1){15}}\put(0,15){\line(1,-1){15}}\put(0,-15){\line(1,1){15}}\put(0,-15){\line(-1,1){15}}\put(15,0){\line(-1,0){30}}\put(17,-2){\tiny$(6)$}
\end{picture}
\begin{picture}(45,20)
\put(15,0){\line(-1,2){5}}\put(15,0){\line(-1,-2){5}}\put(0,15){\line(2,-1){10}}\put(0,15){\line(-2,-1){10}}\put(-15,0){\line(1,2){5}}\put(-15,0){\line(1,-2){5}}\put(0,-15){\line(2,1){10}}\put(0,-15){\line(-2,1){10}}
\put(0,15){\line(0,-1){30}}\put(10,10){\line(0,-1){20}}\put(10,-10){\line(-1,1){20}}\put(-10,-10){\line(0,1){20}}\put(10,-10){\line(-2,5){10}}\put(-10,10){\line(2,-5){10}}\put(17,-2){\tiny$(7)$}
\end{picture}
\begin{picture}(45,20)
\put(15,0){\line(-1,2){5}}\put(15,0){\line(-1,-2){5}}\put(0,15){\line(2,-1){10}}\put(0,15){\line(-2,-1){10}}\put(-15,0){\line(1,2){5}}\put(-15,0){\line(1,-2){5}}\put(0,-15){\line(2,1){10}}\put(0,-15){\line(-2,1){10}}
\put(0,15){\line(0,-1){30}}\put(15,0){\line(-1,1){15}}\put(10,-10){\line(-1,1){20}}\put(-15,0){\line(1,-1){15}}\put(10,-10){\line(-2,5){10}}\put(-10,10){\line(2,-5){10}}\put(17,-2){\tiny$(8)$}
\end{picture}
\begin{picture}(45,20)
\put(15,0){\line(-1,2){5}}\put(15,0){\line(-1,-2){5}}\put(0,15){\line(2,-1){10}}\put(0,15){\line(-2,-1){10}}\put(-15,0){\line(1,2){5}}\put(-15,0){\line(1,-2){5}}\put(0,-15){\line(2,1){10}}\put(0,-15){\line(-2,1){10}}
\put(0,15){\line(0,-1){30}}\put(10,10){\line(-1,-1){20}}\put(15,0){\line(-1,0){30}}\put(15,0){\line(-1,-1){15}}\put(-15,0){\line(1,1){15}}\put(17,-2){\tiny$(9)$}
\end{picture}
\begin{picture}(15,20)
\put(15,0){\line(-1,2){5}}\put(15,0){\line(-1,-2){5}}\put(0,15){\line(2,-1){10}}\put(0,15){\line(-2,-1){10}}\put(-15,0){\line(1,2){5}}\put(-15,0){\line(1,-2){5}}\put(0,-15){\line(2,1){10}}\put(0,-15){\line(-2,1){10}}
\put(0,15){\line(0,-1){30}}\put(10,10){\line(-1,-1){20}}\put(15,0){\line(-1,0){30}}\put(10,-10){\line(-1,1){20}}\put(17,-2){\tiny$(10)$}
\end{picture}&\ \ \ \ \ {\scriptstyle 20+2+4+4+4+1=35={7\choose3}}
\end{array}\]

\vskip1em{\bf\XDCD\ }We identify the set vertices of $R_{l}$ with $\zzz/l\zzz$. Then a map $\alpha:((\zzz/l\zzz)\Delta)_0\to D(R_{l})$ is well-defined by sending $(i,j)$ or $(i,j)_{\pm}$ to the diagonal connecting $i$ and $j$. If $\Delta=D_{m+1}$, then $\alpha$ corresponds $\sigma$-variant points to main diagonals.

Recall that $\tau_2$ gives an automorphism of $(\zzz/l\zzz)\Delta$ of index $2$. We denote by the same letter $\tau_2$ the automorphism of $D(R_{l})$ which is the identity for $\Delta=A_m$ and the $\pi$-radian rotation for other $\Delta$. Then $\alpha$ commutes with $\tau_2$. For any $x\in((\zzz/l\zzz)\Delta)_0$, we can regard $H^-(x)$ and $H^+(x)$ as subsets of $((\zzz/l\zzz)\Delta)_0$. It is not difficult to check the following proposition.

\vskip.5em{\bf Proposition }{\it
(1) $\alpha$ induces a bijection $\alpha:(\zzz\Delta)_0/\langle\tau_2,\sigma\rangle\to D(R_{l})/\langle\tau_2\rangle$.

(2) Fix $x,y\in((\zzz/l\zzz)\Delta)_0$ such that $x\neq y,\tau_2y$ and at most one of $x,y$ is $\sigma$-invariant. Then $y\notin H^-(\tau x)$ if and only if $\alpha(x)$ and $\alpha(y)$ do not cross except their endopoints.}

\vskip.5em{\bf\XDCE\ Theorem }{\it
(1) If $\Delta=A_m$, $B_m$ or $C_m$, then $\alpha$ gives a bijection from the set of maximal $1$-orthogonal subsets of $(\zzz\Delta)_0$ to that of $D(R_{l})$.

(2) If $\Delta=D_{m+1}$, then $\alpha$ gives a surjection from the set of maximal $1$-orthogonal subsets of $(\zzz\Delta)_0$ to that of $D(R_{l})$. For a maximal $1$-orthogonal subset $S$ of $D(R_{l})$, $\#\alpha^{-1}(S)=1$ if $S$ contains only one main diagonal, and $\#\alpha^{-1}(S)=2$ otherwise.}

\vskip.5em{\sc Proof }
(1) By \XDCD(1), $\alpha$ gives a bijection from the set of $\tau_2$-invariant subsets $S$ of $(\zzz\Delta)_0$ to that of $D(R_{l})$. By \XDCD(2), a $\tau_2$-invariant subsets $S$ of $(\zzz\Delta)_0$ satisfies $S\cap(\bigcup_{x\in S}H^-(\tau x))=\emptyset$ if and only if $\alpha(S)$ satisfies \XDCA(i). A maximal $1$-orthogonal subset of $(\zzz\Delta)_0$ is nothing but a maximal $\tau_2$-invariant subset of $(\zzz\Delta)_0$ satisfying $S\cap(\bigcup_{x\in S}H^-(\tau x))=\emptyset$ by \XDBA(5). Obviously, a maximal $1$-orthogonal subset of $D(R_{l})$ is nothing but a maximal $\tau_2$-invariant subset of $D(R_{l})$ satisfying \XDCA(i). Thus the assertion follows.

(2) We denote by $\Sigma$ the set of $\sigma$-variant points of $(\zzz\Delta)_0$. Then $\Sigma$ has two $(\tau\sigma)$-orbits $\Sigma_1$ and $\Sigma_2$. It is not difficult to check that any maximal $1$-orthogonal subset $S$ of $(\zzz\Delta)_0$ satisfies $\#(S\cap\Sigma/\langle\tau_2\rangle)\ge2$ and exactly one of the following conditions.

\strut\kern1em(a) $S\cap\Sigma/\langle\tau_2\rangle=\{x,\sigma x\}$ for some $x$.

\strut\kern1em(b) $S\cap\Sigma\subseteq\Sigma_1$.

\strut\kern1em(c) $S\cap\Sigma\subseteq\Sigma_2$.

By (1) for $\Delta=B_m$ or $C_m$, $\alpha$ gives a bijection from the set of maximal $1$-orthogonal subsets of $(\zzz\Delta)_0$ satisfying (a) above to maximal $1$-orthogonal subsets of $D(R_{l})$ containing only one main diagonal. A similar argument as in the proof of (1) shows that $\alpha$ gives a bijection from the set of maximal $1$-orthogonal subsets of $(\zzz\Delta)_0$ satisfying (b) (resp. (c)) above to the set of maximal $1$-orthogonal subsets of $D(R_{l})$ containing at least two main diagonals. Thus the assertion follows.\rule{5pt}{10pt}

\vskip.5em{\bf\XDD\ }
In the rest of this section, we will give a proof of \XDBB\ and \XDBC. For a set $S$, we denote by $\zzz[S]$ (resp. $\nnn[S]$) the free abelian group (resp. monoid) with the base $S$. For example, the set of isomorphism classes of objects in ${}_\Lambda\dn{M}$ can be identifies with $\nnn[\ind{}_\Lambda\dn{M}]$ by the Krull-Schmidt theorem. For $x=\sum_{y\in S}a_yy\in\zzz[S]$ ($a_y\in\zzz$), put $\supp x:=\{y\in S\ |\ a_y\neq0\}$. We write $x=x_+-x_-$ with $x_+,x_-\in\nnn[S]$ and $\supp x_+\cap\supp x_-=\emptyset$.

Since ${}_\Lambda\underline{\dn{M}}$ forms a $\tau$-category by \XCDE, we can apply the general theory of $\tau$-categories developped in [I1,2,3]. We denote by $\theta X\to X$ the sink map of $X$ in ${}_\Lambda\underline{\dn{M}}$. Define $\theta_i:\nnn[\ind{}_\Lambda\underline{\dn{M}}]\to\nnn[\ind{}_\Lambda\underline{\dn{M}}]$ by $\theta_0X:=X$, $\theta_1X:=\theta X$ and $\theta_iX:=(\theta\theta_{i-1}X-\tau\theta_{i-2}X)_+$ for $i\ge2$. By the theorem [I1;7.1] below, $\theta_i$ becomes a monoid morphism.

\vskip.5em{\bf\XDDA\ Theorem }{\it
For any $X\in{}_\Lambda\underline{\dn{M}}$ and $i\ge0$, the projective cover of the ${}_\Lambda\underline{\dn{M}}$-module $J_{{}_\Lambda\underline{\dn{M}}}^i(\ ,X)$ has the form $\underline{\hom}_\Lambda(\ ,\theta_iX)\to J_{{}_\Lambda\underline{\dn{M}}}^i(\ ,X)\to0$. In particular, $\bigcup_{i\ge0}\supp\theta_iX=\{ Y\in\ind{}_\Lambda\underline{\dn{M}}\ |\ \underline{\hom}_\Lambda(Y,X)\neq0\}$ holds if $\#\ind{}_\Lambda\dn{M}<\infty$.}

\vskip.5em{\bf\XDDB\ }
For any $x\in(\zzz\Delta)_0$, put $\theta x:=\sum_{y\in(\zzz\Delta)_0}a_yy\in\nnn[(\zzz\Delta)_0]$ where $y$ is a predecessor of $x$ with a valued arrow $y\stackrel{(a_y\ *)}{\longrightarrow}x$. We extend $\theta$ to an element of $\endm_{\zzz}(\zzz[(\zzz\Delta)_0])$, and define $\theta_i:\nnn[(\zzz\Delta)_0]\to\nnn[(\zzz\Delta)_0]$ by $\theta_0x:=x$, $\theta_1x:=\theta x$ and $\theta_ix:=(\theta\theta_{i-1}x-\tau\theta_{i-2}x)_+$ for $i\ge2$. An easy calculation gives us the following proposition.

\vskip.5em{\bf Proposition }{\it 
For any $x\in(\zzz\Delta)_0$, $\bigcup_{i\ge0}\supp\theta_ix=H^-(x)$, $\theta_{l-4}x=\tau^{-1}\omega x$, $\theta_{l-3}x=0$ and $\theta_ix=\theta\theta_{i-1}x-\tau\theta_{i-2}x$ for any $i\neq l-2$.}

\vskip.5em{\bf\XDDC\ }In the rest of this paper, we use the notation in \XDBB, and put $G:=\prod_kG_k$. For $X\in\ind{}_\Lambda\dn{M}=\coprod_k(\zzz\Delta_k/G_k)_0$, put $H^{\pm}(X):=\bigcup_{x\in p^{-1}(X)}H^{\pm}(x)$. Thus $H^{\pm}(p(x))=GH^{\pm}(x)$ holds for any $x\in\coprod_k(\zzz\Delta_k)_0$.

\vskip.5em{\bf Proposition }{\it
Let $x,y\in\coprod_k(\zzz\Delta_k)_0$ and $X:=p(x),Y:=p(y)\in\ind{}_\Lambda\dn{M}=\coprod_k(\zzz\Delta_k/G_k)_0$.

(1) $p\tau=\tau p$ and $p\theta_i=\theta_i p$ for any $i\ge0$.

(2) $\underline{\hom}_\Lambda(p(x),p(y))=0$ if and only if $Gx\cap H^-(y)=\emptyset$.

(3) $\underline{\hom}_\Lambda(X,Y)\neq0$ if and only if $p^{-1}(X)\subseteq H^-(Y)$.

(4) For $i>0$, $\ext^i_\Lambda(X,Y)\neq0$ if and only if $p^{-1}(Y)\subseteq H^-(\tau_iX)$.

(5) $p\Omega=\omega p$ and $p\tau_i=\tau_i p$ for any $i>0$.}

\vskip.5em{\sc Proof }
(1) Obviously $p$ commutes with both of $\theta$ and $\tau$. Although $p$ does not commute with $(\ )_+$ in general, \XDDB\ shows that $p$ commutes with $\theta_i$ for any $i$.

(2) Since $p$ commutes with $\theta_i$ by (1), $p(H^-(y))=p(\bigcup_{i\ge0}\supp\theta_iy)=\bigcup_{i\ge0}\supp\theta_iY=\{ Z\in\ind{}_\Lambda\underline{\dn{M}}\ |\ \underline{\hom}_\Lambda(Z,Y)\neq0\}$ holds by \XDDB\ and \XDDA. Thus $Gx\cap H^-(y)=\emptyset$ if and only if $p(x)\notin p(H^-(y))$ if and only if $\underline{\hom}_\Lambda(X,Y)=0$.

(3) Since $(\zzz\Delta)_0$ is a disjoint union of $G$-orbits, the assertion follows from (2).

(4) Since $\ext^i_\Lambda(X,Y)\simeq D\underline{\hom}_\Lambda(Y,\tau_iX)$ holds by \XAE, the assertion follows from (3).

(5) The ${}_\Lambda\underline{\dn{M}}$-module $\underline{\hom}_\Lambda(\ ,X)$ has a simple socle such that $(\soc_{{}_\Lambda\underline{\dn{M}}}\underline{\hom}_\Lambda(\ ,X))(\tau^-\Omega X)\neq0$ by \XBF(1). Since $p$ commutes with $\theta_i$ by (1), $J_{{}_\Lambda\underline{\dn{M}}}^{l-4}(\ ,X)=\soc_{{}_\Lambda\dn{M}}\underline{\hom}_\Lambda(\ ,X)$ and $\tau^-\Omega X=p(\theta_{l-4}x)=p(\tau^{-1}\omega x)$ hold by \XDDA\ and \XDDB. Thus $\Omega X=p(\omega x)$ holds by (1).\rule{5pt}{10pt}


\vskip.5em{\bf\XDDD\ Proof of \XDBB\ and \XDBC\ }We will prove \XDBB. By \XBBB, $\cc$ is maximal $n$-orthogonal if and only if $\ind\underline{\cc}=\bigcap_{X\in\ind\underline{\cc},\ 0<i\le n}\{Y\in\ind{}_\Lambda\underline{\dn{M}}\ |\ \ext^i_\Lambda(X,Y)=0\}$ if and only if $\ind{}_\Lambda\underline{\dn{M}}\backslash\ind\underline{\cc}=\bigcup_{X\in\ind\underline{\cc},\ 0<i\le n}\{Y\in\ind{}_\Lambda\underline{\dn{M}}\ |\ \ext^i_\Lambda(X,Y)\neq0\}$. This is equivalent to $\coprod_k(\zzz\Delta_k)_0\backslash p^{-1}(\ind\underline{\cc})=\bigcup_{X\in\ind\underline{\cc},\ 0<i\le n}H^-(\tau_iX)$ by \XDDC(4). Since $H^-(\tau_iX)=\bigcup_{x\in p^{-1}(\tau_iX)}H^-(x)=\bigcup_{x\in p^{-1}(X)}H^-(\tau_ix)$ holds by \XDDC(5), we obtain \XDBB.

We will prove \XDBC. If $\Delta=A_m$, $B_m$ or $C_m$, then the assertion follows from \XDCB\ and \XDCE(1). Assume $\Delta=D_{m+1}$. The number of maximal $1$-orthogonal subsets of $D(R_{l})$ containing only one main diagonal equals to ${2m\choose m}$ by \XDCB\ for $B_m$, and the number of others equals to ${2m+1\choose m}-{2m\choose m}=\frac{m}{m+1}{2m\choose m}$ by \XDCB\ again. By \XDCE(2), the number of maximal $1$-orthogonal subsets of $(\zzz\Delta)_0$ is ${2m\choose m}+\frac{2m}{m+1}{2m\choose m}=\frac{3m+1}{m+1}{2m\choose m}$.\rule{5pt}{10pt}

\vskip.5em{\bf 5 Appendix }

It is well-known in homological algebra that there exists a bijection between $\ext^n_\Lambda(X,Y)$ and the set of Yoneda classes of exact sequences of length $n$ [HS]. For $n>1$, Yoneda classes are bigger than isomorphism classes of exact sequences. The aim of this section is to show that, for maximal $(n-1)$-orthogonal subcategories, Yoneda classes precisely coincide with isomorphism classes of certain exact sequences.

We keep the notation in \XAD. Fix a maximal $(n-1)$-orthogonal subcategory $\cc$ of ${}_\Lambda\dn{M}$ ($n\ge 1$) and $X,Y\in\cc$. We call an exact sequences ${\bf A}:0\to Y\stackrel{f_n}{\to}C_{n-1}\stackrel{f_{n-1}}{\to}\cdots\stackrel{f_1}{\to}C_0\stackrel{f_0}{\to}X\to0$ with terms in $\cc$ {\it almost-minimal} if $f_i\in J_{\cc}$ holds for any $i$ ($0<i<n$). We say that almost-minimal sequences ${\bf A}$ and ${\bf A}^\prime$ {\it equivalent} if there exists the following commutative diagram whose vertical maps are isomorphisms.
\[\begin{diag}
{\bf A}:\ &0&\RA{}&Y&\RA{f_n}&C_{n-1}&\RA{f_{n-1}}&&\RA{f_2}&C_1&\RA{f_1}&C_0&\RA{f_0}&X&\RA{}&0\\
&&&\parallel&&\downarrow&&\cdots&&\downarrow&&\downarrow&&\parallel\\
{\bf A}^\prime:\ &0&\RA{}&Y&\RA{f_n^\prime}&C_{n-1}^\prime&\RA{f_{n-1}^\prime}&&\RA{f_2^\prime}&C_1^\prime&\RA{f_1^\prime}&C_0^\prime&\RA{f_0^\prime}&X&\RA{}&0
\end{diag}\]

\vskip-.2em{\bf 5.1 Proposition }{\it
Taking Yoneda classes, we have a bijection from the set of equivalence classes of almost-minimal sequences to $\ext^n_\Lambda(X,Y)$.}

\vskip.5em{\sc Proof }
(i) We will show the surjectivity. For $\alpha\in\ext^n_\Lambda(X,Y)$, define a map $\underline{\cc}(\ ,X)\stackrel{\phi}{\to}\ext^n_\Lambda(\ ,Y)$ by $\phi(1_X)=\alpha$. Then $F:=\Im\phi\in\mod\underline{\cc}$ holds by \XCEB. By \XCEA\ and \XCB, there exists an exact sequence ${\bf A}:0\to C_{n+1}\stackrel{f_{n+1}}{\to}\cdots\stackrel{f_1}{\to}C_0\to0$ which induces a minimal projective resolution $0\to\cc(\ ,C_{n+1})\stackrel{\cdot f_{n+1}}{\to}\cdots\stackrel{\cdot f_1}{\to}\cc(\ ,C_0)\to F\to0$ of $F\in\mod\cc$ and a minimal injective resolution $0\to F\to\ext^n_\Lambda(\ ,C_{n+1})$ of $F\in\mod\underline{\cc}$. In particular, we can put $X\simeq X^\prime\oplus C_0$ and $Y\simeq Y^\prime\oplus C_{n+1}$ for some $X^\prime,Y^\prime\in\cc$. Then ${\bf B}:={\bf A}\oplus(0\to Y^\prime\stackrel{1}{\to}Y^\prime\to0\to\cdots\to0\to X^\prime\stackrel{1}{\to}X^\prime\to0)$ is almost-minimal and $\alpha=[{\bf B}]$.

(ii) We will show the injectivity. Let ${\bf A}$ and ${\bf A}^\prime$ be almost-minimal sequences with $[{\bf A}]=[{\bf A}^\prime]$. By \XCB, there exist the following commutative diagrams of exact sequences.
\[\begin{diag}
{\bf A}:\ &0&\RA{}&Y&\RA{f_n}&C_{n-1}&\RA{f_{n-1}}&C_{n-2}&\RA{f_{n-2}}&&\RA{f_1}&C_0&\RA{f_0}&X&\RA{}&0\\
&&&\downarrow^{a_n}&&\downarrow^{a_{n-1}}&&\downarrow^{a_{n-2}}&&\cdots&&\downarrow^{a_0}&&\parallel\\
{\bf A}^\prime:\ &0&\RA{}&Y&\RA{f_n^\prime}&C_{n-1}^\prime&\RA{f_{n-1}^\prime}&C_{n-2}^\prime&\RA{f_{n-2}^\prime}&&\RA{f_1^\prime}&C_0^\prime&\RA{f_0^\prime}&X&\RA{}&0
\end{diag}\]
\[\begin{diag}
\ext^n_\Lambda(X,\ )&\LA{\delta}&\cc(Y,\ )&\LA{f_n\cdot}&\cc(C_{n-1},\ )\\
\parallel&&\uparrow^{a_n\cdot}&&\uparrow^{a_{n-1}\cdot}\\
\ext^n_\Lambda(X,\ )&\LA{\delta^\prime}&\cc(Y,\ )&\LA{f_n^\prime\cdot}&\cc(C_{n-1}^\prime,\ )
\end{diag}\]

Then $\delta(1_Y)=[{\bf A}]=[{\bf A}^\prime]=\delta^\prime(1_Y)$ implies $\delta(a_n-1_Y)=0$. Thus there exists $s\in\cc(C_{n-1},Y)$ such that $a_n=1_Y+f_ns$. Then the following diagram is commutative.
\[\begin{diag}
{\bf A}:\ &0&\RA{}&Y&\RA{f_n}&C_{n-1}&\RA{f_{n-1}}&C_{n-2}&\RA{f_{n-2}}&&\RA{f_1}&C_0&\RA{f_0}&X&\RA{}&0\\
&&&\parallel&&\downarrow^{a_{n-1}-sf_n^\prime}&&\downarrow^{a_{n-2}}&&\cdots&&\downarrow^{a_0}&&\parallel\\
{\bf A}^\prime:\ &0&\RA{}&Y&\RA{f_n^\prime}&C_{n-1}^\prime&\RA{f_{n-1}^\prime}&C_{n-2}^\prime&\RA{f_{n-2}^\prime}&&\RA{f_1^\prime}&C_0^\prime&\RA{f_0^\prime}&X&\RA{}&0
\end{diag}\]

Using almost-minimality, we can easily check that vertical maps in the diagram above are isomorphisms. Thus ${\bf A}$ and ${\bf A}^\prime$ are equivalent.\rule{5pt}{10pt}

\vskip.5em{\footnotesize
\begin{center}
{\bf References}
\end{center}

[ArS] M. Artin, W. F. Schelter: Graded algebras of global dimension $3$. Adv. in Math. 66 (1987), no. 2, 171--216.

[AV] M. Artin, J.-L. Verdier: Reflexive modules over rational double points. Math. Ann. 270 (1985), no. 1, 79--82.

[A1] M. Auslander: Coherent functors. 1966 Proc. Conf. Categorical Algebra (La Jolla, Calif., 1965) pp. 189--231 Springer, New York.


[A2] M. Auslander: Functors and morphisms determined by objects. Representation theory of algebras (Proc. Conf., Temple Univ., Philadelphia, Pa., 1976), pp. 1--244. Lecture Notes in Pure Appl. Math., Vol. 37, Dekker, New York, 1978.

[A3] M. Auslander: Isolated singularities and existence of almost split sequences. Representation theory, II (Ottawa, Ont., 1984), 194--242, Lecture Notes in Math., 1178, Springer, Berlin, 1986.

[A4] M. Auslander: Rational singularities and almost split sequences. Trans. Amer. Math. Soc. 293 (1986), no. 2, 511--531. 

[AB] M. Auslander, M. Bridger: Stable module theory. Memoirs of the American Mathematical Society, No. 94 American Mathematical Society, Providence, R.I. 1969.

[ABu] M. Auslander, R. Buchweitz: The homological theory of maximal Cohen-Macaulay approximations. Colloque en l'honneur de Pierre Samuel (Orsay, 1987). Mem. Soc. Math. France (N.S.) No. 38, (1989), 5--37.

[AR1] M. Auslander, I. Reiten: Stable equivalence of dualizing $R$-varieties. Advances in Math. 12 (1974), 306--366.

[AR2] M. Auslander, I. Reiten: Almost split sequences in dimension two. Adv. in Math. 66 (1987), no. 1, 88--118.


[ARS] M. Auslander, I. Reiten, S. O. Smalo: Representation theory of Artin algebras. Cambridge Studies in Advanced Mathematics, 36. Cambridge University Press, Cambridge, 1995. 

[AS] M. Auslander, S. O. Smalo: Almost split sequences in subcategories. J. Algebra 69 (1981), no. 2, 426--454. 

[CE] H. Cartan, S. Eilenberg: Homological algebra. Princeton University Press, Princeton, N. J., 1956.

[CR] C. W. Curtis, I. Reiner: Methods of representation theory. Vol. I. With applications to finite groups and orders. A Wiley-Interscience Publication. John Wiley \& Sons, Inc., New York, 1990.

[DKR] Y. A. Drozd, V. V. Kiri\v cenko, A. V. Ro\u\i ter: Hereditary and Bass orders. (Russian) Izv. Akad. Nauk SSSR Ser. Mat. 31 (1967), 1415--1436.

[E] H. Esnault: Reflexive modules on quotient surface singularities. J. Reine Angew. Math. 362 (1985), 63--71. 

[EG] E. G. Evans, P. Griffith: Syzygies. London Mathematical Society Lecture Note Series, 106. Cambridge University Press, Cambridge, 1985. 

[FGR] R. M. Fossum, P. A. Griffith, I. Reiten: Trivial extensions of abelian categories. Homological algebra of trivial extensions of abelian categories with applications to ring theory. Lecture Notes in Mathematics, Vol. 456. Springer-Verlag, Berlin-New York, 1975. 

[GR] P. Gabriel, A. V. Roiter: Representations of finite-dimensional algebras. Springer-Verlag, Berlin, 1997. 

[GL] W. Geigle, H. Lenzing: A class of weighted projective curves arising in representation theory of finite-dimensional algebras. Singularities, representation of algebras, and vector bundles (Lambrecht, 1985), 265--297, Lecture Notes in Math., 1273, Springer, Berlin, 1987. 

[GN1] S. Goto, K. Nishida: Finite modules of finite injective dimension over a Noetherian algebra. J. London Math. Soc. (2) 63 (2001), no. 2, 319--335. 

[GN2] S. Goto, K. Nishida: Towards a theory of Bass numbers with application to Gorenstein algebras. Colloq. Math. 91 (2002), no. 2, 191--253. 

[H] D. Happel: Triangulated categories in the representation theory of finite-dimensional algebras. London Mathematical Society Lecture Note Series, 119. Cambridge University Press, Cambridge, 1988.

[HPR] D. Happel, U. Preiser, C. M. Ringel: Vinberg's characterization of Dynkin diagrams using subadditive functions with application to $D{\rm Tr}$-periodic modules. Representation theory, II (Proc. Second Internat. Conf., Carleton Univ., Ottawa, Ont., 1979), pp. 280--294, Lecture Notes in Math., 832, Springer, Berlin, 1980.

[HS] P. J. Hilton, U. Stammbach: A course in homological algebra. Graduate Texts in Mathematics, 4. Springer-Verlag, New York, 1997.

[I1] O. Iyama: $\tau$-categories I: Ladders, Algebr. Represent. Theory 8 (2005), no. 3, 297--321.

[I2] O. Iyama: $\tau$-categories II: Nakayama pairs and Rejective subcategories, Algebr. Represent. Theory 8 (2005), no. 4, 449--477.

[I3] O. Iyama: $\tau$-categories III: Auslander orders and Auslander-Reiten quivers, Algebr. Represent. Theory 8 (2005), no. 5, 601--619.


[Mc] J. McKay: Graphs, singularities, and finite groups. The Santa Cruz Conference on Finite Groups (Univ. California, Santa Cruz, Calif., 1979), pp. 183--186, Proc. Sympos. Pure Math., 37, Amer. Math. Soc., Providence, R.I., 1980. 

[M] H. Matsumura: Commutative ring theory. Cambridge Studies in Advanced Mathematics, 8. Cambridge University Press, Cambridge, 1989.


[R1] C. Riedtmann: Representation-finite self-injective algebras of class $A\sb{n}$. Representation theory, II (Proc. Second Internat. Conf., Carleton Univ., Ottawa, Ont., 1979), pp. 449--520, Lecture Notes in Math., 832, Springer, Berlin, 1980. 

[R2] C. Riedtmann: Algebren, Darstellungskocher, Uberlagerungen und zuruck. (German) Comment. Math. Helv. 55 (1980), no. 2, 199--224.

[Ri] C. M. Ringel: Tame algebras and integral quadratic forms. Lecture Notes in Mathematics, 1099. Springer-Verlag, Berlin, 1984. 

[RV] I. Reiten, M. Van den Bergh: Two-dimensional tame and maximal orders of finite representation type. Mem. Amer. Math. Soc. 80 (1989).

[RS] K. W. Roggenkamp, J. W. Schmidt: Almost split sequences for integral group rings and orders. Comm. Algebra 4 (1976), no. 10, 893--917.

[S] D. Simson: Linear representations of partially ordered sets and vector space categories. Algebra, Logic and Applications, 4. Gordon and Breach Science Publishers, Montreux, 1992. 

[St] R. P. Stanley: Enumerative combinatorics. Vol. 2. Cambridge Studies in Advanced Mathematics, 62. Cambridge University Press, Cambridge, 1999. 


[V] M. Van den Bergh: Non-commutative crepant resolutions. The legacy of Niels Henrik Abel, 749--770, Springer, Berlin, 2004.

[Y] Y. Yoshino: Cohen-Macaulay modules over Cohen-Macaulay rings. London Mathematical Society Lecture Note Series, 146. Cambridge University Press, Cambridge, 1990. 

[W1] A. Wiedemann: Classification of the Auslander-Reiten quivers of local Gorenstein orders and a characterization of the simple curve singularities. J. Pure Appl. Algebra 41 (1986), no. 2-3, 305--329.

[W2] A. Wiedemann: Die Auslander-Reiten Kocher der gitterendlichen Gorensteinordnungen. Bayreuth. Math. Schr. No. 23 (1987), 1--134. 
}

\vskip.5em{\footnotesize
{\sc Department of Mathematics, University of Hyogo, Himeji, 671-2201, Japan}

{\it iyama@sci.u-hyogo.ac.jp}


\vskip.5em
Current address: 

{\sc Graduate School of Mathematics, Nagoya University,

Chikusa-ku, Nagoya, 464-8602, Japan}

{\it iyama@math.nagoya-u.ac.jp}}
\end{document}